\documentclass[11pt]{article}
\usepackage{amsfonts}
\usepackage{amsfonts,mathrsfs,amssymb,amsthm,mathptm}
\usepackage{amsmath,amscd}
\usepackage{mathptm,pslatex}
\usepackage{color}
\usepackage[dvips]{graphicx}
\usepackage[all]{xy}
\usepackage{graphicx}
\usepackage{fancyhdr}
\begin{document}

\newtheorem{Def}{Definition}[section]
\newtheorem{Bsp}[Def]{Example}
\newtheorem{Prop}[Def]{Proposition}
\newtheorem{Theo}[Def]{Theorem}
\newtheorem{Lem}[Def]{Lemma}
\newtheorem{Koro}[Def]{Corollary}
\theoremstyle{definition}
\newtheorem{Rem}[Def]{Remark}

\newcommand{\add}{{\rm add}}
\newcommand{\con}{{\rm con}}

\newcommand{\colim}{{\rm colim\, }}
\newcommand{\limt}{{\rm lim\, }}
\newcommand{\Add}{{\rm Add }}
\newcommand{\Tor}{{\rm Tor}}
\newcommand{\Cogen}{{\rm Cogen}}
\newcommand{\Tria}{{\rm Tria}}
\newcommand{\tria}{{\rm tria}}
\newcommand{\E}{{\rm E}}
\newcommand{\Mor}{{\rm Morph}}
\newcommand{\End}{{\rm End}}
\newcommand{\ind}{{\rm ind}}
\newcommand{\rsd}{{\rm res.dim}}
\newcommand{\rd} {{\rm rep.dim}}
\newcommand{\ol}{\overline}
\newcommand{\overpr}{$\hfill\square$}
\newcommand{\rad}{{\rm rad}}
\newcommand{\soc}{{\rm soc}}
\renewcommand{\top}{{\rm top}}
\newcommand{\pd}{{\rm proj.dim}}
\newcommand{\id}{{\rm inj.dim}}
\newcommand{\fld}{{\rm flat.dim}}
\newcommand{\Fac}{{\rm Fac}}
\newcommand{\Gen}{{\rm Gen}}
\newcommand{\gd}{{\rm gl.dim}}
\newcommand{\dm}{{\rm dom.dim}}
\newcommand{\tdim}{{\rm dim}}
\newcommand{\fd} {{\rm fin.dim}}
\newcommand{\Fd} {{\rm Fin.dim}}
\newcommand{\Pf}[1]{{\mathscr P}^{<\infty}(#1)}
\newcommand{\DTr}{{\rm DTr}}
\newcommand{\cpx}[1]{#1^{\bullet}}
\newcommand{\D}[1]{{\mathscr D}(#1)}
\newcommand{\Dz}[1]{{\mathscr D}^+(#1)}
\newcommand{\Df}[1]{{\mathscr D}^-(#1)}
\newcommand{\Db}[1]{{\mathscr D}^b(#1)}
\newcommand{\C}[1]{{\mathscr C}(#1)}
\newcommand{\Cz}[1]{{\mathscr C}^+(#1)}
\newcommand{\Cf}[1]{{\mathscr C}^-(#1)}
\newcommand{\Cb}[1]{{\mathscr C}^b(#1)}
\newcommand{\Dc}[1]{{\mathscr D}^c(#1)}
\newcommand{\K}[1]{{\mathscr K}(#1)}
\newcommand{\Kz}[1]{{\mathscr K}^+(#1)}
\newcommand{\Kf}[1]{{\mathscr  K}^-(#1)}
\newcommand{\Kb}[1]{{\mathscr K}^b(#1)}
\newcommand{\modcat}{\ensuremath{\mbox{{\rm -mod}}}}
\newcommand{\Modcat}{\ensuremath{\mbox{{\rm -Mod}}}}

\newcommand{\stmodcat}[1]{#1\mbox{{\rm -{\underline{mod}}}}}
\newcommand{\pmodcat}[1]{#1\mbox{{\rm -proj}}}
\newcommand{\imodcat}[1]{#1\mbox{{\rm -inj}}}
\newcommand{\Pmodcat}[1]{#1\mbox{{\rm -Proj}}}
\newcommand{\Imodcat}[1]{#1\mbox{{\rm -Inj}}}
\newcommand{\opp}{^{\rm op}}
\newcommand{\otimesL}{\otimes^{\rm\mathbb L}}
\newcommand{\rHom}{{\rm\mathbb R}{\rm Hom}\,}
\newcommand{\projdim}{\pd}
\newcommand{\flatdim}{\fld}
\newcommand{\Hom}{{\rm Hom}}
\newcommand{\Coker}{{\rm Cok}}
\newcommand{ \Ker  }{{\rm Ker}}
\newcommand{ \Cone }{{\rm Con}}
\newcommand{ \Img  }{{\rm Im}}
\newcommand{\Ext}{{\rm Ext}}
\newcommand{\StHom}{{\rm \underline{Hom}}}

\newcommand{\gm}{{\rm _{\Gamma_M}}}
\newcommand{\gmr}{{\rm _{\Gamma_M^R}}}

\def\vez{\varepsilon}\def\bz{\bigoplus}  \def\sz {\oplus}
\def\epa{\xrightarrow} \def\inja{\hookrightarrow}

\newcommand{\lra}{\longrightarrow}
\newcommand{\lraf}[1]{\stackrel{#1}{\lra}}
\newcommand{\ra}{\rightarrow}
\newcommand{\dk}{{\rm dim_{_{k}}}}

{\Large \bf
\begin{center}
Recollements of derived categories III: finitistic dimensions

\end{center}}

\medskip
\centerline{\textbf{Hong Xing Chen} and \textbf{Chang Chang Xi}$^*$}

\renewcommand{\thefootnote}{\alph{footnote}}
\setcounter{footnote}{-1} \footnote{ $^*$ Corresponding author.
Email: xicc@cnu.edu.cn; Fax: +86 10 68903637.}
\renewcommand{\thefootnote}{\alph{footnote}}
\setcounter{footnote}{-1} 
\renewcommand{\thefootnote}{\alph{footnote}}
\setcounter{footnote}{-1} 

\begin{abstract}
In this paper, we study homological dimensions of algebras linked by recollements of derived module categories, and establish a series of new upper bounds and relationships among their finitistic or global dimensions. This is closely related to a longstanding conjecture, the finitistic dimension conjecture, in representation theory and homological algebra. Further, we apply our results to a series of situations of particular interest: exact contexts, ring extensions, trivial extensions, pullbacks of rings, and algebras induced from Auslander-Reiten sequences. In particular, we not only extend and amplify Happel's reduction techniques for finitistic dimenson conjecture to more general contexts, but also generalise some recent results in the literature.

\end{abstract}

\tableofcontents

\section{Introduction}
Recollements of triangulated categories have been introduced by Beilinson, Bernstein and Deligne in order to decompose derived categories of sheaves into two parts, an open and a closed one (see \cite{BBD}), and thus providing a natural habitat for Grothendieck's six functors. Similarly, recollements of derived module categories can be seen as short exact sequences, describing a derived module category in terms of a subcategory and of a quotient, both of which may be derived module categories themselves, related by six functors that in general are not known. It turns out that recollements provide a very useful framework
for understanding connections among three algebraic or geometric objects in which one is interested.

In a series of papers on recollements of derived module categories, we are addressing basic questions about recollements and the rings involved.
Our starting point has been infinite-dimensional tilting theory (see \cite{xc1}). While Happel's theorem establishes a derived equivalence between a given ring and the endomorphism ring of a
finitely generated tilting object (see \cite{Happel1, CPS}), Bazzoni has shown that for large tilting modules one gets instead a recollement relating three triangulated categories, with
two of them being the derived categories of the given ring and the endomorphism ring of the large tilting module. In \cite{xc1} we have addressed the question of determining
the third category in this recollement as a derived category of a ring and we have explained this ring in terms of universal localisations in the sense of Cohn (see \cite{cohnbook1, nr} for definition). Among the applications has been a counterexample to the
Jordan-H\"older problem for derived module categories. In \cite{xc3} we have dealt with the problem of constructing recollements in order to relate
rings. Our main construction, of exact contexts, can be seen as a far-reaching generalisation of pullbacks of rings. In \cite{xc5} we have used this
construction to relate algebraic $K$-theory of different rings. It turned out that under mild assumptions,
the $K$-theory of an algebra can be fully decomposed under a sequence of recollements.

For cohomology and for homological invariants of algebras, such a complete decomposition is not possible. Nevertheless,
results by Happel in \cite{Happel}
for the case of bounded derived categories (when fewer recollements exist than in the unbounded case) show that finiteness of  finitistic dimension
of an algebra can be reduced along a recollement; if such an invariant is finite for the two outer terms, then it is finite for the middle term, too. Note that
the particular values of these invariants depend on the ring and are not invariants of the derived category. The present paper aims at extending Happel's
reduction techniques for homological conjectures.  As in Happel's paper \cite{Happel} we will focus on finitistic dimensions, which include finite global dimensions as a special case.

Recall that the finitistic dimension
of a ring $R$, denoted by fin.dim$(R)$, is by definition the supremum of projective dimensions of those left $R$-modules having a finite projective resolution by finitely generated projective modules. The well-known finitistic dimension conjecture states that any
Artin algebra should have finite finitistic dimension (see, for instance, \cite[conjecture (11), p.410]{ars}). This conjecture is a longstanding question (\cite[Bass, 1960]{bass}) and has still not been settled. It is closely related to at least seven other main conjectures in homological representation theory of algebras (see \cite[p. 409-410]{ars}). In the literature, there is another definition of big finitistic dimension of a ring $R$, denoted by Fin.dim$(R)$, which is the supremum of projective dimensions of all those left $R$-modules which have finite projective dimension. Clearly, fin.dim$(R)\le$ Fin.dim$(R)$. Usually, they are quite different (see \cite{ZH}).

There are two main directions in this article. First, we provide reduction techniques for homological invariants of unbounded derived module categories, that is, for the most general possible setup (which also has been covered in the preceding articles in this series). In the first main result, Theorem \ref{main-1}, we give criteria for the finiteness of finitistic dimension for each of the three rings in a recollement of derived module categories, in terms of the other two. The criteria aim to be applicable by putting conditions on particular objects, not on the whole category. The second main result, Theorem \ref{hom-dim}, applies the first
main result to the general contexts of the so-called exact contexts introduced in \cite{xc1}, and in addition provides upper and lower bounds for the finitistic dimensions of the three rings involved. A series of corollaries then applies the general results to classes of examples of particular interest, such as ring extensions, trivial extensions, quotient rings and endomorphism rings of modules related by an almost split sequence.


\begin{Theo}\label{main-1}
Let $R_1$, $R_2$ and $R_3$ be rings. Suppose that there exists a recollement among
the derived module categories $\D{R_3}$, $\D{R_2}$ and $\D{R_1}$ of $R_3$, $R_2$ and $R_1:$

$$\xymatrix@C=1.2cm{\D{R_1}\ar[r]^-{{i_*}}
&\D{R_2}\ar[r]^-{\;j^!}\ar@/^1.2pc/[l]^-{i^!}\ar_-{i^*}@/_1.2pc/[l]
&\D{R_3} \ar@/^1.2pc/[l]^{j_*\;}\ar@/_1.2pc/[l]_{j_!\;}}$$
Then the following hold true:

$(1)$ Suppose that $j_!$ restricts to a functor $\Db{R_3}\to \Db{R_2}$ of bounded derived module categories. If $\fd(R_2)<\infty$, then $\fd(R_3)<\infty$.

$(2)$ Suppose that  $i_*(R_1)$ is a compact object in $\D {R_2}$. Then we have the following:

\;\quad $(a)$ If $\fd(R_2)<\infty$, then $\fd(R_1)<\infty$.

\;\quad $(b)$ If $\fd(R_1)<\infty$ and $\fd(R_3)<\infty,$ then $\fd(R_2)<\infty$.
\end{Theo}

Note that the assumption of Theorem \ref{main-1} on unbounded derived module
categories is weaker than the one on bounded derived module categories, because
the existence of recollements of bounded derived module categories implies the one of unbounded derived module categories.
This is shown by a recent investigation on recollements at different levels in \cite{AKLY, koenig}.
So, Theorem \ref{main-1} (see also Corollary \ref{db-level}) generalizes the main result
in \cite{Happel} since for a recollement of $\Db{R_j\modcat}$ with $R_j$ a finite-dimensional
algebra over a field for $1\le j\le 3$, one can always deduce that $i_*(R_1)$ is compact in $\D{R_2}$.
Moreover, Theorem \ref{main-1} extends and amplifies a result in \cite{DMX}
because we deal with arbitrary rings instead of Artin algebras, and also yields a generalization of a result in \cite{XP} for left coherent rings to the one for arbitrary rings (see Corollary \ref{f3-0} below).

To prove this result, we introduce homological widths (or cowidth)
for complexes that are quasi-isomorphic to bounded complexes of
projective (or injective) modules (see Section \ref{sect3.1} for
details). Broadly speaking, the homological width (respectively,
cowidth) defines a map from homotopy equivalence classes of bounded
complexes of projective (respectively, injective) modules to the
natural numbers. It measures, up to homotopy equivalence, how large
the minimal interval of such a complex is in which its non-zero
terms are distributed. Particularly, if a module has finite
projective dimension, then its homological width is exactly the
projective dimension. Using homological widths, we will present a
substantial and technical result, Theorem \ref{finitistic}, which is a strengthened version of Theorem \ref{main-1} and
describes explicitly upper bounds for finitistic dimensions, so that
Theorem \ref{main-1} will become an easy consequence of Theorem
\ref{finitistic}. Note that, in \cite{Happel}, one of the key arguments in proofs is that finite-dimensional algebras have finitely many non-isomorphic simple modules, while in our general context we do not have this fact and therefore must avoid this kind of arguments. So, the idea of proving Theorems \ref{finitistic} and \ref{main-1} will be completely
different from the ones in \cite{Happel} and
\cite{DMX}. Moreover, our methods also lead to results on upper bounds for big finitistic  and global dimensions. For details, we refer the reader to Theorems \ref{big finitistic} and \ref{gldim}.

\medskip
Now, let us utilize Theorem \ref{main-1} to recollements constructed in \cite{xc3} and establish relationships among finitistic dimensions of noncommutative tensor products and related rings. First of all, we  recall some notions from \cite{xc3}:

Let $R$, $S$ and $T$ be associative rings with identity, and let $\lambda:R\to S$ and $\mu:R\to T$ be ring homomorphisms.
Suppose that $M$ is an $S$-$T$-bimodule together with an element $m\in M$. We say that the quadruple $(\lambda, \mu, M, m)$ is an \emph{exact context} if the following sequence
$$
0\lra R\lraf{(\lambda,\,\mu)}S\oplus T
\lraf{\left({\cdot\,m\,\atop{-m\,\cdot}}\right)}M\lra 0
$$
is an exact sequence of abelian groups, where $\cdot m$ and $m \cdot$ denote the right and left multiplication by $m$ maps, respectively. There is a list of examples in \cite{xc3} that guarantees the ubiquity of exact contexts.

Given an exact context $(\lambda,\mu, M, m)$, there is defined a ring with identity in \cite{xc3}, called the \emph{noncommutative tensor product} of $(\lambda,\mu, M, m)$ and denoted by  $T\boxtimes_RS$ if the meaning of the exact context is clear. This notion not only generalizes the one of usual tensor products over commutative rings and captures coproducts of rings, but also plays a key role in describing the left parts of recollements induced from homological exact contexts (see \cite[Theorem 1.1]{xc3}).

For an $R$-module $_RX$, we denote by $\fld(_RX)$ and $\pd(_RX)$ the flat and projective dimensions of $X$, respectively.

From the proof of Theorem \ref{main-1} or Theorem \ref{finitistic}, we have the following result.

\begin{Theo} \label{hom-dim}
Let $(\lambda, \mu, M, m)$ be an exact context with the noncommutative tensor product $T\boxtimes_RS$. Then

$(1)$ $\fd(R)\leq \fd(S)+\fd(T)+\max\{1,\fld(T_R)\}+1.$

$(2)$ Suppose that $\Tor^R_i(T,S)=0$ for all $i\geq 1$. If the left $R$-module $_RS$ has a finite projective resolution by finitely generated projective modules, then the following hold
true:

\quad $(a)$ $\fd(T\boxtimes_RS)\leq \fd(S)+\fd(T)+1$.

\quad $(b)$ $\fd(S)\leq \fd\left(\begin{array}{lc} S& M\\ 0 & T\end{array}\right)\leq \fd(R)+\fd(T\boxtimes_RS)+\max\{1,\pd(_RS)\}+3$.
\end{Theo}

Note that for the triangular matrix algebra $B:=\left(\begin{array}{lc} S& M\\ 0 & T\end{array}\right)$, it is known that $\fd(B)\le \fd(S)+\fd(T)+1$.
But, Theorem \ref{hom-dim}(2)(b) provides us with a new upper bound for the finitistic dimension of $B$. That is, the finiteness of $\fd(B)$ can be seen from the one of $\fd(T\boxtimes_RS)$ and $\fd(R)$, involving the starting ring $R$ but without information on $\fd(S)$ and $\fd(T)$. This is non-trivial and somewhat surprising. Moreover, in Theorem \ref{hom-dim}(2), if $\lambda:R\ra S$ is a homological ring epimorphism, then we even obtain better estimations: $\fd(S)\leq \fd(R)$ and $\fd(T\boxtimes_RS)\leq\fd(T)$. In this case, $T\boxtimes_RS$ can be interpreted as the coproduct $S\sqcup_RT$ of the $R$-rings of $S$ and $T$.

\smallskip
Now, let us state several consequences of Theorem \ref{hom-dim}. First, we utilize Theorem \ref{hom-dim} to finitistic dimensions of ring extensions. This is of particular interest because the finitistic dimension conjecture can be reformulated over perfect fields in terms of ring extensions (see \cite{xx}). Note that, in the following result, we do not impose any conditions on the radicals of rings, comparing with \cite{x2, xx}.

\begin{Koro}\label{ring extension}
Suppose that $S\subseteq R$ is an extension of rings, that is, $S$ is a subring of $R$ with the same identity. Let $R'$ be the endomorphism ring of the $S$-module $R/S$, and let $R'\boxtimes_SR$ be the noncommutative tensor product of the exact context determined the extension. Then

 $(1)$  $\fd(S)\leq \fd(R)+\fd(R')+\max\big\{1,\,\fld((R/S)_S),\,\fld\big(\Hom_S(R,R/S)_S\big)\big\}+1.$

 $(2)$ Suppose that the left $S$-module $R$ is projective and finitely generated.  Then the following hold true:

\quad  $(a)$ $\fd(R'\boxtimes_SR)\leq \fd(R)+\fd(R')+1$.

\quad  $(b)$ $\fd(R)\leq \fd(S)+\fd(R'\boxtimes_SR)+4$.
\end{Koro}

Next, we apply Theorem \ref{hom-dim} to trivial extensions. Recall that, given a ring $R$ and an $R$-$R$-bimodule $M$, the \emph{trivial extension} of $R$ by $M$ is a ring, denoted by $R\ltimes M$, with abelian group $R\oplus M$ and multiplication: $(r,m)(r',m')=(rr',rm'+mr')$ for $r,r'\in R$ and $m,m'\in M$. For consideration of $\Fd(R\ltimes M)$, we refer the reader to \cite[Chapter 4]{FGR}.

\begin{Koro}\label{mod1b}
Let $\lambda: R\to S$ be a ring epimorphism and $M$ an $S$-$S$-bimodule such that $\Tor^R_i(M,S)$ = $0$ for all $i\geq 1$. If $_RS$ has a finite projective resolution by finitely generated projective $R$-modules, then

\quad $(a)$ $\fd(S\ltimes M)\leq \fd(S)+\fd(R\ltimes M)+1$.

\quad $(b)$ $\fd(S)\leq \fd(R)+\fd(S\ltimes M)$.

\end{Koro}

Now, we apply Theorem \ref{hom-dim} to pullback squares of rings and surjective homomorphisms.

\begin{Koro}\label{mod1a}
Let $R$ be a ring, and let $I_1$ and $I_2$ be ideals of $R$ such that $I_1\cap I_2=0$. Then

$(1)$ $\fd(R)\leq \fd(R/I_1)+\fd(R/I_2)+\max\{1,\fld((R/I_2)_R)\}+1.$

$(2)$ Suppose that $\Tor_i^R(I_2,I_1)=0$ for all $i\geq 0$. If the left $R$-module $R/I_1$ has a finite projective resolution by finitely generated projective modules, then

\quad $(a)$ $\fd(R/(I_1+I_2))\leq \fd(R/I_1)+\fd(R/I_2)+1$.

\quad $(b)$ $\fd(R/I_1)\leq \fd(R)+\fd(R/(I_1+I_2))+\max\{1,\pd(_R(R/I_1))\}+3$.
\end{Koro}

The strategy of proving Corollaries \ref{mod1b} and \ref{mod1a} is as follows: First, we show that under the given assumptions we can get exact pairs, a class of special exact contexts, and then employ Theorem \ref{hom-dim} by verifying the Tor-vanishing condition. At last, we have to describe  noncommutative tensor products more substantially for the cases considered.

Finally, we mention a corollary on finitistic dimensions of algebras arising from idempotent ideals and almost split sequences (see \cite{ars} for definition).

\begin{Koro} \label{ars-dim}
$(1)$ If $I$ is an idempotent ideal in a ring $R$, then  $\fd(R/I)\le \fd\big(\End_R(R\oplus I)\big)\leq \fd\big(\End_R(_RI)\big)+\fd(R/I)+2.$

$(2)$  Let $0\ra Z\ra Y\ra X\ra 0$ be an almost split sequence of $R$-modules with $R$ an Artin algebra. If $\Hom_R(Y,Z)=0$, then
$\fd(\End_R(Y\oplus X))\leq \fd(\End_R(Y))+2.$
\end{Koro}

The paper is sketched as follows: In Section \ref{sect2}, we first recall some necessary definitions and then prove two results on coproducts of rings. In Section \ref{sect3}, we provide all proofs of our results. Especially, we introduce homological widths of complexes and prove an amplified version, Theorem \ref{finitistic}, of Theorem \ref{main-1} phrased in terms of homological widths and finitistic dimensions of involved rings, such that Theorem \ref{main-1} is deduced readily from Theorem \ref{finitistic}. Moreover, the methods developed in this section also give similar upper bounds for global and big finitistic dimensions (see Theorems \ref{big finitistic} and \ref{gldim}).

\section{Definitions and conventions\label{sect2}}
In this section, we fix notation and briefly recall some definitions. For unexplained ones, we refer the reader to \cite{xc3, xc5}.

Throughout the paper, all notation and terminology are standard. For example, by a ring we mean an associative ring with identity. For a ring $R$, we denote by $R$-Mod the category of all left $R$-modues, and by $\C{R}$, $\K {R}$ and $\D R$  the unbounded complex, homotopy and derived categories of $R$-Mod, respectively. As usual, by adding a superscript $*\in \{-, +, b\}$, we denote their corresponding $*$-bounded categories, for instance, $\Db R$ is the bounded derived category of $R$-Mod. The full subcategory of compact objects in $\D{R}$ is denoted by $\Dc{R}$. This category
is also called the \emph{perfect derived module category} of $R$. It is known that the localization functor $\K{R}\to\D{R}$ induces a triangle equivalence from the homotopy category of bounded complexes of finitely generated projective $R$-modules to $\Dc{R}$.

As usual, we write a complex in $\C{R}$ as $\cpx{X}=(X^i, d_{\cpx X}^i)_{i\in\mathbb{Z}}$, and call $d_{\cpx X}^i:X^i\to X^{i+1}$ the $i$-th differential of $\cpx{X}$. Sometimes, for simplicity, we shall write $(X^i)_{i\in\mathbb{Z}}$ for $\cpx{X}$ without mentioning the morphisms $d^i_{\cpx X}$. Given a chain map $\cpx{f}:\cpx{X}\to\cpx{Y}$ in $\C{R}$, its mapping cone is denoted by $\Cone(\cpx{f})$. For an integer $n$, the $n$-th cohomology of $\cpx{X}$ is denoted by $H^n(\cpx{X})$. Let $\sup{(\cpx{X})}$ and $\inf{(\cpx{X})}$ be the supremum and minimum of indices $i\in \mathbb{Z}$ such that $H^i(\cpx{X})\ne 0$, respectively. If $\cpx{X}$ is acyclic, that is, $H^i(\cpx{X})=0$ for all $i\in\mathbb{Z}$, then we understand that $\sup(\cpx{X})=-\infty$ and $\inf(\cpx{X})=+\infty$. If $\cpx{X}$ is not acyclic, then  $\inf(\cpx{X})\leq \sup(\cpx{X})$, and $H^n(\cpx{X})=0$ if $\sup{(\cpx{X})}$ is an integer and $n>\sup{(\cpx{X})}$ or if $\inf{(\cpx{X})}$ is an integer and $n<\inf{(\cpx{X})}$.
For a complex $\cpx{X}$ in $\C{R}$, if $H^n(\cpx{X})=0$ for almost all $n$, then $\cpx{X}$ is isomorphic in $\D{R}$ to a bounded complex. So,  $\Db{R}$ is equivalent to the full subcategory of $\D{R}$ consisting of all complexes with finitely many nonzero cohomologies.

As a convention, we write the composite of two homomorphisms $f:X\ra Y$ and $g: Y\ra Z$ in $R$-Mod as $fg$. Thus the image of an element $x\in X$ under $f$ will be written on the opposite of the scalars as $(x)f$ instead of $f(x)$. This convention makes $\Hom_R(X,Y)$ naturally a left $\End_R(X)$- and right $\End_R(Y)$-bimodule. But, for two functors $F: \mathcal{X}\ra \mathcal{Y}$ and $G:\mathcal{Y}\ra \mathcal{Z}$ of categories, we write $GF:\mathcal{X}\ra \mathcal{Z}$ for their composition.

Let us now recall the notion of recollements of triangulated
categories, which was defined by Beilinson, Bernstein and Deligne
in \cite{BBD} to study derived categories of
perverse sheaves over
singular spaces. It may be thought as a kind of categorifications of exact sequences in abelian categories.

\begin{Def}\label{def01} \rm
Let  $\mathcal{D}$, $\mathcal{D'}$ and $\mathcal{D''}$ be
triangulated categories. We say that $\mathcal{D}$ is a
\emph{recollement} of $\mathcal{D'}$ and $\mathcal{D''}$ if there
are six triangle functors among the three categories:
$$\xymatrix{\mathcal{D''}\ar^-{i_*=i_!}[r]&\mathcal{D}\ar^-{j^!=j^*}[r]
\ar^-{i^!}@/^1.2pc/[l]\ar_-{i^*}@/_1.6pc/[l]
&\mathcal{D'}\ar^-{j_*}@/^1.2pc/[l]\ar_-{j_!}@/_1.6pc/[l]}$$ such
that

$(1)$ $(i^*,i_*),(i_!,i^!),(j_!,j^!)$ and $(j^*,j_*)$ are adjoint
pairs,

$(2)$ $i_*,j_*$ and $j_!$ are fully faithful functors,

$(3)$ $j^! i_!=0$ (and thus also $i^!j_*=0$ and $i^*j_!=0$), and

$(4)$ for each object $X\in\mathcal{D}$, there are two triangles in
$\mathcal D$ induced by counit and unit adjunctions:
$$
i_!i^!(X)\lra X\lra j_*j^*(X)\lra i_!i^!(X)[1],$$ $$ j_!j^!(X)\lra
X\lra i_*i^*(X)\lra j_!j^!(X)[1],
$$where the shift functor of triangulated categories is denoted by [$1$].
\end{Def}

Recall that the \emph{coproduct} of a family $\{R_i\mid i\in I\}$ of $R_0$-rings
with $I$ an index set is defined to be an $R_0$-ring $R$ together
with a family $\{\rho_i: R_i\ra R\mid i\in I \}$ of
$R_0$-homomorphisms of rings such that, for  any $R_0$-ring $S$ with a family
of $R_0$-homomorphisms $\{\tau_i: R_i\ra S\mid i\in I\}$, there
exists a unique $R_0$-homomorphism $\delta: R\ra S$ such that
$\tau_i=\rho_i\delta$ for all $i\in I$. It is well known that coproducts of rings exist (see \cite{cohn}). However, this existence result does not provide us with a handy form of coproducts; therefore we need a concrete description of coproducts for our situations considered.

In the following we describe coproducts of rings for two cases in terms of some known constructions. This will be used in later proofs. The first one is for trivial extensions

\begin{Lem} \label{tri-ext}
Suppose that $\lambda: R\to S$ is a ring epimorphism and $M$ is an
$S$-$S$-bimodule. Let $\widetilde{\lambda}: R\ltimes M\to S\ltimes
M$ be the ring homomorphism induced by
$\lambda$. Then the coproduct $S\sqcup_{R}(R\ltimes M)$ is
isomorphic to $S\ltimes M$, that is, the inclusion $S\ra S\ltimes M$
and $\widetilde{\lambda}$ define the coproduct.
\end{Lem}

{\it Proof.} Let $\mu: R\to R\ltimes M$  and $\rho: S\to S\ltimes M$
be the inclusions of rings. Note that $S$ and $R\ltimes M$ are
$R$-rings via $\lambda$ and $\mu$, respectively, and that
$\lambda\rho=\mu\,\widetilde{\lambda}: R\to S\ltimes M$. To prove
that $S\ltimes M$, together with $\rho$ and $\widetilde{\lambda}$,
is the coproduct of $S$ and $R\ltimes M$ over $R$, we suppose that $\Lambda$ is an arbitrary ring and that $f:
R\ltimes M\to \Lambda$ and $g: S\to \Lambda$ are arbitrary ring
homomorphisms such that $\lambda\,g=\mu\,f$. Then we have to show that there is a unique ring homomorphism $h: S\ltimes
M\to \Lambda$ such that $\widetilde{\lambda}\,h=f$ and $\rho\, h=g$.
Clearly, if such an $h$ exists, then $h$ must be defined by $(s,
m)\mapsto (m)f +(s)g$ for $s\in S$ and $m\in M$. This shows the
uniqueness of $h$. So, it suffices to show that the above-defined map
$h$ is a ring homomorphism. Certainly, $h$ is a homomorphism of
abelian groups. Hence, we have to show that $h$ preserves multiplication.

Let $s_i\in S$ and $m_i\in M$ for $i=1, 2$. On the one hand,
$\big((s_1, m_1)\, (s_2, m_2)\big)h=(s_1s_2,\, s_1m_2+ m_1s_2)h=(
s_1m_2+ m_1s_2)f+(s_1s_2)g=(s_1m_2) f+ (m_1s_2)f+(s_1)g \,(s_2)g.$
On the other hand, $\big((s_1, m_1)\big)h\,\big((s_2,
m_2)\big)h=\big((m_1)f +(s_1)g\big)\big((m_2)f
+(s_2)g\big)=(m_1)f\,(m_2)f + (m_1)f\,(s_2)g +(s_1)g\,(m_2)f +
(s_1)g \,(s_2)g=(m_1m_2)f+ (m_1)f (s_2)g +(s_1)g (m_2)f + (s_1)g
\,(s_2)g= (m_1)f (s_2)g +(s_1)g\, (m_2)f + (s_1)g\, (s_2)g$ since
$m_1m_2=0$. This implies that if $(s_1m_2)f=(s_1)g\, (m_2)f$ and
$(m_1s_2)f=(m_1)f\, (s_2)g $, then $\big((s_1, m_1)\, (s_2,
m_2)\big) h$ =$((s_1, m_1))h $ $((s_2, m_2))h.$ So, to prove that $h$
preserves multiplication, we need only to verify these additional
conditions, that is,
$$(sm)f=(s)g \,(m)f \; \mbox{ and } \; (ms)f=(m)f\,(s)g \, \mbox{ for }
\, s\in S \, \mbox{  and } \, m\in M.$$

To show the former, we fix an $m\in M$ and define two maps:
$$\varphi: S\to \Lambda,\;  s\mapsto (sm)f \;\,\mbox{and}\;\,
\psi: S\to \Lambda,\; s\mapsto (s)g\,(m)f.$$  Since $\lambda g=\mu\,f$, one can check that both
$\varphi$ and $\psi$ are homomorphisms of $R$-modules such that
$\lambda \varphi=\lambda\psi$. But we do not know if they are homomorphisms of rings. Nevertheless, we can still have $\phi=\psi$ because $\lambda: R\to S$ being a ring epimorphism by assumption implies that  the map $\Hom_R(R,\lambda): \Hom_R(S,\Lambda)\ra \Hom_R(R,\Lambda)$ is an isomorphism, and therefore it is injective. Thus $\phi=\psi$. Similarly, we can show that
$(ms)f=(m)f\,(s)g$. Consequently, the map $h$ preserves
multiplication and is actually a ring homomorphism.
$\square$

\medskip
The other description of coproducts is for quotients of rings by ideals, which applies to Milnor squares (see \cite{Milnor}).

\begin{Lem}\label{epi00'}
Let $R_0$ be a ring, and let $R_i$ be an $R_0$-ring  with ring
homomorphism $\lambda_i: R_0\ra R_i$ for $i=1,2$.

$(1)$ If $\,\lambda_1: R_0\to R_1$ is a ring epimorphism, then so
is the canonical homomorphism $\rho_2: R_2\ra R_1\sqcup_{R_0}R_2$.

$(2)$ Let $I$ be an ideal of $R_0$, and let $J$ be the ideal of
$R_2$ generated by the image  $(I)\lambda_2$ of $I$ under the
map $\lambda_2$. If $R_1=R_0/I$ and $\lambda_1: R_0\to R_1$
is the canonical surjective map, then $R_1\sqcup_{R_0}R_2=R_2/J$.
\end{Lem}

{\it Proof.} $(1)$ It follows from the definition of coproducts of
rings that $\lambda_1\,\rho_1=\lambda_2\rho_2:R_0\to
R_1\sqcup_{R_0}R_2$. We point out that $\rho_2$ is a ring
epimorphism. In fact, if $f,g: R_1\sqcup_{R_0}R_2\ra S$ are two ring
homomorphisms such that $\rho_2f=\rho_2g$, then
$\lambda_2\rho_2f=\lambda_2\rho_2g$. This means that
$\lambda_1\rho_1 f=\lambda_1\rho_1 g$, and therefore $\rho_1
f=\rho_1 g$ since $\lambda_1$ is a ring epimorphism. By the
universal property of coproducts, we have $g=f$. Thus $\rho_2$ is a
ring epimorphism.

$(2)$ Let $\rho_2:R_2\to R_2/J$ be the canonical surjection, and let
$\rho_1:R_1\to R_2/J$ be the ring homomorphism induced by
$\lambda_2$ since
$J=R_2\,(I)\lambda_2\,R_2\supseteq(I)\lambda_2$. Now, we
claim that $R_2/J$ together with $\rho_1$ and $\rho_2$ is the
coproduct of $R_1$ and $R_2$ over $R_0$. Clearly, we have
$\lambda_1\,\rho_1=\lambda_2\rho_2:R_0\to R_2/J$. Further,
assume that $\tau_1: R_1\ra S$ and $\tau_2:R_2\to S$ are two ring
homomorphisms such that $\lambda_2\tau_2=\lambda_1\tau_1$.
Then $(I)\lambda_2\tau_2=(I)\lambda_1\tau_1=0$, and
therefore $(J)\tau_2=0$. This means that there is a unique ring
homomorphism $\delta:R_2/J\to S$ such that $\tau_2=\rho_2\delta$. It
follows that $\lambda_1\tau_1=\lambda_2\tau_2=
\lambda_2\rho_2\delta=\lambda_1\rho_1\delta$. Since
$\lambda_1$ is surjective, we have $\tau_1=\rho_1\delta$. This
shows that $R_1\sqcup_{R_0}R_2=R_2/J$.  $\square$

\section{Proofs \label{sect3}}
This section is devoted to proofs of all results mentioned in the introduction. We start with introducing the so-called homological widths for complexes, and then prove a strengthened version, Theorem \ref{finitistic} below, of Theorem \ref{main-1}. As consequences of Theorem \ref{finitistic}, we get proofs of Theorems \ref{main-1} and \ref{hom-dim}. Finally, we apply Theorem \ref{hom-dim} to give proofs of all corollaries, and mention two results on global and big finitistic dimensions.

\subsection{Homological widths and cowidths of complexes}\label{sect3.1}
As a generalization of finite projective or injective dimensions of modules, we define, in this subsection, homological widths and cowidths for bounded complexes of projective and injective modules, respectively.

Let $R$ be a ring. For an $R$-module $M$, we denote by $\pd(M)$, $\id(M)$ and $\fld(M)$ the projective, injective and flat dimension of $M$, respectively. As usual, $\Pmodcat{R}$ is the category of all projective left $R$-modules, and $\pmodcat{R}$ is the full subcategory of $\Pmodcat{R}$ consisting of all finitely generated projective left $R$-modules. If there is a projective resolution  $0\ra P_n\ra \cdots \ra P_1\ra P_0\ra M\ra 0$ of $M$ with all $P_i$ in $\pmodcat{R}$, then we say that $M$ is of \emph{finite type}. The category of all $R$-modules of finite type will be denoted by $\Pf{R}$.

Let $\cpx{P}:=(P^n,d^n_{\cpx P})_{n\in\mathbb{Z}}\in\Cb{\Pmodcat{R}}$. We define the homological \emph{width} of $\cpx{P}$ in the following way:
$$
w(\cpx{P}):=\left\{\begin{array}{ll} 0 & \mbox{if}\;\cpx{P}\;\; \mbox{is acyclic},\\ \sup(\cpx{P})-\inf(\cpx{P})+\pd\big(\Coker(d_{\cpx P}^{\,\inf(\cpx{P})-1})\big)
& \mbox{otherwise. }\end{array} \right.
$$
Clearly, $0\leq w(\cpx{P})<\infty$. Moreover, $\cpx{P}$ is isomorphic in $\Kb{\Pmodcat{R}}$ to a complex
$$
\cpx{Q}:\; 0\lra Q^{t-p}\lraf{d^{t-p}} Q^{t-p+1}\lraf{d^{t-p+1}} \cdots \lra Q^{t-1}\lraf{d^{t-1}} Q^{t} \lra P^{t+1}\lraf{d_{\cpx{P}}^{t+1}} \cdots\lra P^{s-1}\lraf{d_{\cpx{P}}^{s-1}} \Ker(d_{\cpx{P}}^s)\lra 0
$$
with $s:=\sup(\cpx{P})$, $t:=\inf(\cpx{P})$,
$p:=\pd\big(\Coker(d_{\cpx P}^{t-1})\big)$ and each term being projective. Clearly, the sequence
$$
0\lra Q^{t-p}\lraf{d^{t-p}} Q^{t-p+1}\lraf{d^{t-p+1}} \cdots \lra Q^{t-1}\lraf{d^{t-1}} Q^{t}\lra \Coker(d_{\cpx P}^{t-1})\lra 0
$$
is a projective resolution of the $R$-module $\Coker(d_{\cpx P}^{t-1})$. Note that if $\cpx{P}\in\Cb{\pmodcat{R}}$, we can choose $\cpx{Q}\in\Cb{\pmodcat{R}}$.

\smallskip
The following result says that homological widths of bounded complexes of projective modules
are preserved under homotopy equivalences.

\begin{Lem}\label{equ}
Let $\cpx{M},\cpx{N}\in\Cb{\Pmodcat{R}}$. If $\cpx{M}\simeq \cpx{N}$ in $\Kb{\Pmodcat{R}}$, then
$w(\cpx{M})=w(\cpx{N})$.
\end{Lem}

{\it Proof.} Recall that $\Kb{\Pmodcat R}$ is the stable category of the Frobenius category $\Cb{\Pmodcat R}$ with projective objects being acyclic complexes. Assume that $\cpx{M}\simeq \cpx{N}$ in $\Kb{\Pmodcat{R}}$. Then there exist two acyclic complexes $\cpx{P}$ and $\cpx{Q}$ in $\Cb{\Pmodcat R}$ such that $\cpx{M}\oplus \cpx{P}\simeq \cpx{N}\oplus \cpx{Q}$ in $\Cb{\Pmodcat
R}$. This implies that $H^i(\cpx{M})\simeq H^i(\cpx{N})$ and $\Coker(d_{\cpx M}^i)\oplus \Coker(d_{\cpx P}^i)\simeq \Coker(d_{\cpx N}^i)\oplus \Coker(d_{\cpx Q}^i)$ for all $i\in\mathbb{Z}$. Thus $\sup(\cpx{M})=\sup(\cpx{N})$ and $\inf(\cpx{M})=\inf(\cpx{N})$. Moreover, since $\Coker(d_{\cpx P}^i)$ and $\Coker(d_{\cpx Q}^i)$ belong to $\Pmodcat{R}$, we have $\pd\big(\Coker(d_{\cpx M}^i)\big)=\pd\big(\Coker(d_{\cpx N}^i)\big)$. It follows that $w(\cpx{M})=w(\cpx{N})$. $\square$

\medskip
Thanks to Lemma \ref{equ}, the definition of homological widths for complexes can be extended slightly to derived categories in the following sense: Given a complex $\cpx{X}\in\D{R}$, if there is a complex $\cpx{P}\in\Cb{\Pmodcat{R}}$ such that  $\cpx{X}\simeq \cpx{P}$ in $\D{R}$, then we define $w(\cpx{X}):=w(\cpx{P})$. This is well defined: If there exists another complex $\cpx{Q}\in\Cb{\Pmodcat{R}}$ such that $\cpx{X}\simeq \cpx{Q}$ in $\D{R}$, then $\cpx{P}\simeq \cpx{Q}$ in $\Kb{\Pmodcat{R}}$ and $w(\cpx{P})=w(\cpx{Q})$ by Lemma \ref{equ}. So, for such a complex $\cpx{X}$, its homological width $w(\cpx{X})$ can be characterized as follows:
$$w(\cpx{X})=\min\left\{
\alpha_{\cpx{P}}-\beta_{\cpx{P}}\,\geq 0\, \bigg |
\begin{array}{ll}
\cpx{P}\simeq \cpx{X} \;\;\mbox{in}\;\;\D{R}\;\;\mbox{for}\;\;
\cpx{P}\in\Cb{\Pmodcat{R}} \\
\;\;\mbox{with}\;\; P^i=0\;\;\mbox{for}\;\; i<\beta_{\cpx{P}}\;\;\mbox{or}\;\; i>\alpha_{\cpx{P}}
\end{array}
\right\}.
$$
Clearly, if $X\in R\Modcat$ has finite projective dimension, then $w(X)=\pd(X)$.

Dually, we can define homological cowidths for bounded complexes of injective $R$-modules.

Let $\Imodcat{R}$ denote the category of injective $R$-modules.
Given a complex $\cpx{I}:=(I^n,d^n_{\cpx I})_{n\in\mathbb{Z}}\in\Cb{\Imodcat{R}}$,
we define the homological \emph{cowidth} of $\cpx{I}$ as follows:
$$
cw(\cpx{I}):=\left\{\begin{array}{ll} 0 & \mbox{if}\;\cpx{I}\;\; \mbox{is acyclic},\\ \sup(\cpx{I})-\inf(\cpx{I})+\id\big(\Ker(d_{\cpx I}^{\,\sup(\cpx{I})})\big)
& \mbox{otherwise. }\end{array} \right.
$$
Similarly, if a complex $\cpx{Y}$ is isomorphic in $\D{R}$ to a bounded complex
$\cpx{I}\in\Cb{\Imodcat{R}}$, then we define $cw(\cpx{Y}):=cw(\cpx{I})$.
In particular, if $Y\in R\Modcat$ has finite injective dimension, then $cw(Y)=\id(Y)$.
Also, we have the following characterization of $cw(\cpx{Y})$:
$$
cw(\cpx{Y})=\min\left\{
\alpha_{\cpx{I}}-\beta_{\cpx{I}}\,\geq 0\, \bigg  |
\begin{array}{ll}
\cpx{I}\simeq \cpx{Y} \;\;\mbox{in}\;\;\D{R}\;\;\mbox{for}\;\;
\cpx{I}\in\Cb{\Imodcat{R}} \\
\;\;\mbox{with}\;\; I^i=0\;\;\mbox{for}\;\; i<\beta_{\cpx{I}}\;\;\mbox{or}\;\; i>\alpha_{\cpx{I}}
\end{array}
\right\}.
$$

Homological widths and cowidths will be used to bound homological dimensions in the next section.

\subsection{Proof of Theoem \ref{main-1}}
In this subsection, we shall first prove an amplified version of Theorem \ref{main-1}, namely Theorem \ref{finitistic} below, so that Theorem \ref{main-1} becomes a straightforward consequence of Theorem \ref{finitistic}.

Recall that the \emph{finitistic dimension} of a ring $R$, denoted by
$\fd(R)$, is defined as follows:
$$
\fd(R):= \mbox{sup}{\{\pd(_{R}X)\mid X\in \Pf{R}\}}.
$$

For each  $n\in\mathbb{Z}$, we define
$$
\mathscr{D}^c_{\geq n}(R):=\{\cpx{X}\in\Dc{R}\mid \cpx{X}\simeq
\cpx{P}\;\,\mbox{in}\;\,\Dc{R}\;\,\mbox{with}\;\, \cpx{P}\in
\Cb{\pmodcat {R}}\;\,\mbox{such that}\;\,{P^i}=0\;\,\mbox{for all }\;\,i<n \}.
$$
From this definition, we have $\mathscr{D}^c_{\geq n}(R)\subseteq \mathscr{D}^c_{\geq n'}(R)$ whenever $n\geq n'$. Since the localization functor $\K{R}\to\D{R}$ induces a triangle equivalence $\Kb{\pmodcat{R}}\lraf{\simeq}\Dc{R}$, we have
$$\Dc{R}= \bigcup_{n\in\mathbb{Z}}\mathscr{D}^c_{\geq n}(R).$$ Clearly,
if $\fd(R)=m<\infty$, then $\Pf{R}\subseteq
\mathscr{D}^c_{\geq -m}(R)$. For the convenience of later discussions, we also formally set
$\mathscr{D}^c_{\geq -\infty}(R):=\Dc{R}$ and $\mathscr{D}^c_{\geq +\infty}(R):=\{0\}.$

\begin{Lem} \label{preparation}
Let $m, n\in\mathbb{N}$. Then the following statements are true:

$(1)$ The full subcategory $\mathscr{D}^c_{\geq n}(R)$ of $\Dc{R}$ is closed under direct summands in $\Dc{R}$.

$(2)$ Let $\cpx{X}\in\mathscr{D}^c_{\geq n}(R)$, $\cpx{Z}\in\mathscr{D}^c_{\geq m}(R)$ and $s=\min{\{n,m\}}$. Then, for any distinguished triangle $\cpx{X}\to \cpx{Y}\to\cpx{Z}\to\cpx{X}[1]$ in $\Dc{R}$, we have $\cpx{Y}\in\mathscr{D}^c_{\geq s}(R)$.
 \end{Lem}

{\it Proof.} $(1)$ Let $\cpx{M}\in\Kb{\pmodcat{R}}$, and let $\cpx{N}:=(N^i)_{i\in\mathbb{Z}}\in\Kb{\pmodcat{R}}$ such that $N^i=0$ for all $i<n$.
Suppose that $\cpx{M}$ is a direct summand of $\cpx{N}$ in $\Kb{\pmodcat{R}}$, or equivalently,  there is a complex  $\cpx{L}\in\Cb{\pmodcat{R}}$ such that $\cpx{M}\oplus\cpx{L}\simeq \cpx{N}$ in $\Kb{\pmodcat{R}}$. Hence $H^i(\cpx{M})=0$ for all $i<n$. Note that $\Kb{\pmodcat R}$ is the stable category of the Frobenius category $\Cb{\pmodcat R}$ with projective objects being acyclic complexes. So we can find two acyclic complexes $\cpx{U}$ and $\cpx{V}$ in $\Cb{\pmodcat R}$ such that $\cpx{M}\oplus\cpx{L}\oplus\cpx{U}\simeq \cpx{N}\oplus \cpx{V}$ in $\Cb{\pmodcat R}$. This implies that $$\Coker(d_{\cpx M}^{n-1})\oplus \Coker(d_{\cpx L}^{n-1})\oplus \Coker(d_{\cpx U}^{n-1})\simeq \Coker(d_{\cpx N}^{n-1})\oplus \Coker(d_{\cpx V}^{n-1})=N^n\oplus \Coker(d_{\cpx V}^{n-1}).$$ Since $N^n\oplus \Coker(d_{\cpx V}^{n-1})\in\pmodcat{R}$, we have $\Coker(d_{\cpx M}^{n-1})\in\pmodcat{R}$. It follows that $\cpx{M}$ is isomorphic in
$\Kb{\pmodcat R}$ to the following
truncated complex $$0\lra \Coker(d_{\cpx M}^{n-1})\lra M^{n+1}\lra M^{n+2}
\lra\cdots \lra 0.$$ Recall that the
localization functor $\K{R}\to\D{R}$ induces a triangle equivalence
$\Kb{\pmodcat R}\lraf{\simeq} \Dc{R}$. Thus $(1)$ follows.

$(2)$ Since $\cpx{X}\in\mathscr{D}^c_{\geq n}(R)$, there exists a complex
$\cpx{P}\in\Cb{\pmodcat{R}}$ with $P^i=0$ for $i<n$ such that $\cpx{X}\simeq \cpx{P}$ in $\Dc{R}$. Similarly, there exists another complex $\cpx{Q}\in\Cb{\pmodcat{R}}$ with $Q^i=0$ for $i<m$ such that $\cpx{Z}\simeq \cpx{Q}$ in $\Dc{R}$. It follows from the triangle equivalence  $\Kb{\pmodcat R}\lraf{\simeq} \Dc{R}$ that
$$
\Hom_{\Kb{\pmodcat R}}(\cpx{Q}[-1],\cpx{P})\simeq
\Hom_{\Dc{R}}(\cpx{Q}[-1],\cpx{P})\simeq\Hom_{\Dc{R}}(\cpx{Z}[-1],\cpx{X}).
$$
Thus the given triangle yields a distinguished triangle in $\Dc{R}$:
$$
\cpx{Q}[-1]\lraf{\cpx{f}} \cpx{P}\lra \cpx{Y}\lra \cpx{Q}
$$
with $\cpx{f}$ a chain map in $\C{R}$. Then $\cpx{Y}\simeq \Cone(\cpx{f})$ in $\Dc{R}$. Since $\Cone(\cpx{f})^i=Q^i\oplus P^i$ for any $i\in\mathbb{Z}$, we have $\Cone(\cpx{f})\in\Cb{\pmodcat R}$ and $\Cone(\cpx{f})^i=0$ for $i<s$. This implies $\cpx{Y}\in\mathscr{D}^c_{\geq s}(R)$. $\square$

\medskip
To investigate relationships among finitistic dimensions of rings in recollements, it may be convenient to introduce the notion of \emph{finitistic dimensions of functors}.

Let $R_1$ and $R_2$ be two arbitrary rings. Suppose that $\mathscr{X}_1$ and $\mathscr{X}_2$ are full subcategories of $\D{R_1}$ and $\D{R_2}$, respectively, and that  $R_1\Modcat\subseteq \mathscr{X}_1$. For a given additive functor $F:\mathscr{X}_1\ra \mathscr{X}_2$, we define
$$
\inf{(F)}:=\inf\{n\in\mathbb{Z}\mid H^n(F(X))\neq 0\;\;\mbox{for some}\;\; X\in R_1\Modcat\},
$$
$$
\fd(F):=\inf\{n\in\mathbb{Z}\mid H^n(F(X))\neq 0\;\;\mbox{for some}\;\; X\in\Pf{R_1}\}.
$$
Note that $\inf(F)=+\infty$ if and only if $F(X)=0$ in $\D{R_2}$ for all $X\in R_1\Modcat$. In fact, if there exits some $X\in R_1\Modcat$ such that $H^n(F(X))\neq 0$ for some integer $n$, then $\inf(F)\leq n$. Moreover, by definition, we always have $\inf(F)\leq \fd(F)$ and $\fd(F)\in \mathbb{Z}\cup \{-\infty, +\infty\}$.
\smallskip

\begin{Lem}\label{f0}
Let $F:\D{R_1}\to \D{R_2}$ be a triangle functor. Then the following statements are true:

$(1)$ If $F$ has a left adjoint $L:\D{R_2}\to\D{R_1}$ with $L(R_2)\in\Df{R_1}$, then $\inf(F)\geq -\sup{(L(R_2))}$.

$(2)$ If $F$ has a right adjoint $G:\D{R_2}\to\D{R_1}$, then $F$ restricts to a functor $\Db{R_1}\to \Db{R_2}$ if and only if $G\big(\Hom_\mathbb{Z}(R_2,\mathbb{Q}/\mathbb{Z})\big)$
is isomorphic in $\D{R_1}$ to a bounded complex $\cpx{I}$ of injective $R_1$-modules. In this case, $\inf(F)\geq -\big(m+\id(\Ker(d_{\cpx{I}}^m))\big)$, where $m:=\sup(\cpx{I})$ and $d_{\cpx{I}}^m:I^m\to I^{m+1}$ is the $m$-th differential of $\cpx{I}$.

\end{Lem}

{\it Proof.}
$(1)$ For each $n\in\mathbb{Z}$ and $M\in R_1\Modcat$, we have $$H^n(F(M))\simeq \Hom_{\D{R_2}}(R_2, F(M)[n])\simeq\Hom_{\D{R_1}}(L(R_2), M[n]).$$
Since $L(R_2)\in\Df{R_1}$, we have $s:=\sup(L(R_2))<+\infty$. Recall that the localization functor $\K{R_1}\to \D{R_1}$ induces a triangle equivalence $\Kf{\Pmodcat{R_1}}\lraf{\simeq} \Df{R_1}$. So there is a complex $\cpx{P}:=(P^j)_{j\in\mathbb{Z}}\in\Cf{\Pmodcat{R_1}}$ with $P^j=0$ for all $j>s$ such that $\cpx{P}\simeq L(R_2)$ in $\D{R_1}$. It follows that $$H^n(F(M))\simeq \Hom_{\D{R_1}}(L(R_2), M[n])\simeq \Hom_{\D{R_1}}(\cpx{P}, M[n])\simeq\Hom_{\K{R_1}}(\cpx{P}, M[n]) =0$$ for all $n<-s$. Thus $\inf(F)\geq -s$.

$(2)$ To calculate cohomologies of complexes, we consider the functor
$$(-)^\vee:=\Hom_\mathbb{Z}(-, \mathbb{Q}/\mathbb{Z}):\;
\mathbb{Z}\Modcat\lra \mathbb{Z}\Modcat.$$
This is an exact functor with the property  that a $\mathbb{Z}$-module $U$ is zero if and only if so is $U^\vee$, because $\mathbb{Q}/\mathbb{Z}$ is an injective cogenerator for $\mathbb{Z}\Modcat$.

Let $\cpx{X}\in\D{R_2}$. Then
{\small $$
{H^0(\cpx{X})}^\vee=\Hom_\mathbb{Z}(H^0(\cpx{X}),\mathbb{Q}/\mathbb{Z})\simeq
\Hom_{\K{\mathbb{Z}}}(\cpx{X},\mathbb{Q}/\mathbb{Z})\simeq\Hom_{\K{\mathbb{Z}}}
(R_2\otimes_{R_2}\cpx{X},\mathbb{Q}/\mathbb{Z})\simeq\Hom_{\K{R_2}}(\cpx{X}, {R_2}^\vee).
$$}
Since ${R_2}^\vee$ is an injective $R_2$-module, we have $\Hom_{\K{R_2}}(\cpx{X}, {R_2}^\vee)\simeq \Hom_{\D{R_2}}(\cpx{X}, {R_2}^\vee)$. Thus $$ {H^0(\cpx{X})}^\vee\simeq\Hom_{\D{R_2}}(\cpx{X}, {R_2}^\vee).$$
Now, let $M\in R_1\Modcat$ and $n\in\mathbb{Z}$. Then $H^n(F(M))^\vee\simeq \Hom_{\D{R_2}}\big(F(M)[n], {R_2}^\vee)$.  Since $(F, G)$ is an adjoint pair, we have $$\Hom_{\D{R_2}}\big(F(M)[n], {R_2}^\vee)\simeq\Hom_{\D{R_1}}(M[n], G({R_2}^\vee)\big).$$
This implies that $H^n(F(M))=0$ if and only if $\Hom_{\D{R_1}}(M[n], G({R_2}^\vee)\big)=0$.

Let $W=G({R_2}^\vee)$. To check the sufficiency of $(2)$, it is enough to show that $\Hom_{\D{R_1}}(M[n], \cpx{W}\big)=0$ for almost all $n$. In fact, if $\cpx{W}$ is isomorphic in $\D{R_1}$ to a bounded complex $\cpx{I}$ of injective $R_1$-modules, then $$\Hom_{\D{R_1}}(M[n], \cpx{W}\big)\simeq \Hom_{\D{R_1}}(M[n], \cpx{I}\big)\simeq \Hom_{\K{R_1}}(M[n], \cpx{I}\big)=0$$
for almost all $n$.

In the following, we will show the necessity of $(2)$. Suppose that $F$ restricts to a functor $\Db{R_1}\to \Db{R_2}$. We first claim that $H^n(\cpx{W})=0$ for almost all $n$, that is $\cpx{W}\in\Db{R_1}$.

Actually, we have the following isomorphisms of abelian groups:
$$
H^n(\cpx{W})\simeq\Hom_{\D{R_1}}(R_1,G(R_2^\vee)[n]) \simeq\Hom_{\D{R_2}}(F(R_1), R_2^\vee[n])\simeq H^{-n}(F(R_1)).
$$
Since $F(R_1)\in\Db{R_2}$, we have $H^n(F(R_1))=0$ for almost all $n$. Thus $H^n(\cpx{W})=0$ for almost all $n$. In other words, $\cpx{W}$ is isomorphic in $\Db{R_1}$ to a bounded complex. Consequently, there exists a lower-bounded complex $\cpx{I}$ of injective $R_1$-modules such that $\cpx{I}\simeq \cpx{W}$ in $\D{R_1}$. In particular, we have $H^n(\cpx{I})\simeq H^n(\cpx{W})$ for all $n$.
To complete the proof of the necessity of $(2)$, it remains to show that $\cpx{I}$ can be chosen to be a bounded complex.

Note that we have the following isomorphisms:
$$
\Hom_{\D{R_2}}(F(M), R_2^\vee[n])\simeq \Hom_{\D{R_1}}(M,\cpx{W}[n])\simeq\Hom_{\D{R_1}}(M,\cpx{I}[n])\simeq\Hom_{\K{R_1}}(M,\cpx{I}[n]).
$$
As $F:\D{R_1}\to\D{R_2}$ restricts to a functor $\Db{R_1}\to \Db{R_2}$ by assumption, we get $F(M)\in\Db{R_2}$. Up to isomorphism in $\D{R_2}$, we may assume that $F(M)\in\Cb{R_2}$. Since $R_2^\vee$ is an injective $R_2$-module, we see that $\Hom_{\D{R_2}}(F(M), R_2^\vee[n])\simeq \Hom_{\K{R_2}}(F(M), R_2^\vee[n])=0$ for almost all $n$.
Thus $\Hom_{\K{R_1}}(M,\cpx{I}[n])=0$ for almost all $n$. Particularly, there is a natural number $\delta_M$ (depending on $M$) such that $\Hom_{\K{R_1}}(M$, $\cpx{I}[n])=0$ for all $n >\delta_M$. We may suppose that the complex $\cpx{I}$ is of the following form:
$$
0\lra I^{s}\lraf{d^s}I^{s+1}\lraf{d^{s+1}} \cdots \lra I^m\lraf{d^{m}} I^{m+1}\lraf{d^{m+1}}\cdots
\lra I^{i}\lraf{d^{i}}I^{i+1}\lra\cdots
$$
where all terms $I^i$ are injective and where $s\leq m:=\sup(\cpx{I})$ and $H^i(\cpx{I})=0$ for any $i>m$. Let $V:=\bigoplus_{i\geq m} \Img(d^i)$. Then $$\Hom_{\K{R_1}}(V,\cpx{I}[n])=0\;\;\mbox{for all}\;\;n>\delta_V.$$
Now we define $t:=\max{\{m,\delta_V\}}$. Then $\Hom_{\K{R_1}}(\Img(d^t),\cpx{I}[t+1])=0$. This implies that the chain map $\Img(d^t)\to \cpx{I}[t+1]$, induced from the inclusion $\Img(d^t)\hookrightarrow I^{t+1},$ is homotopic to the zero map. Therefore, the canonical surjection $I^t\twoheadrightarrow\Img(d^t)$ must split. With $I^t$ then also $\Img(d^t)$ is an injective module. Since $H^i(\cpx{I})=0$ for any $i>m$, we see that $\cpx{I}$ is isomorphic in $\D{R_1}$ to the following bounded complex:
$$
0\lra I^{s}\lraf{d^s}I^{s+1}\lraf{d^{s+1}} \cdots \lra I^m\lraf{d^{m}} I^{m+1}\lraf{d^{m+1}}\cdots
\lra I^{t}\lraf{d^{t}}\Img(d^t)\lra 0
$$
with all of its terms being injective. Thus, up to isomorphism in $\D {R_1}$, we can choose $\cpx{I}$ to be a bounded complex of injective modules. This completes the proof of the necessity of $(2)$.

To show the last statement of $(2)$, we note that the $R_1$-module $\Ker(d^m)$ has a finite injective resolution since $H^i(\cpx{I})=0$ for all $i>m$. Hence, up to isomorphism in $\D{R_1}$, we can replace $\cpx{I}$ by the following bounded complex of injective $R_1$-modules:
$$
0\lra I^{s}\lraf{d^s}I^{s+1}\lraf{d^{s+1}} \cdots \lra I^{m-1}\lra  \widetilde{I}^m\lraf{\widetilde{d}^{m}} \widetilde{I}^{m+1}\lraf{\widetilde{d}^{m+1}}\cdots
\lra \widetilde{I}^{m+p-1}\lraf{\widetilde{d}^{m+p-1}}\widetilde{I}^{m+p}\lra 0
$$ where $\Ker(\widetilde{d}^{m})=\Ker(d^m)$ and $p:=\id(\Ker(d^m))\leq t$.
This implies that $\Hom_{\K{R_1}}(M[n],\cpx{I})=0$ for all $n<-(m+p)$. Since $$
H^n(F(M))^\vee\simeq \Hom_{\D{R_1}}(M[n], \cpx{W})\simeq \Hom_{\D{R_1}}(M[n], \cpx{I})\simeq \Hom_{\K{R_1}}(M[n],\cpx{I}),
$$
we have $H^n(F(M))=0$ for all $n<-(m+p)$. Thus $\inf(F)\geq -(m+p)$. $\square$

\medskip
We remark that, in Lemma \ref{f0}(2), the $R_2$-module $I:=\Hom_{\mathbb{Z}}(R_2,\mathbb{Q}/\mathbb{Z})$ can be replaced by any injective cogenerator of
$R_2\Modcat$. This is due to the fact that $G$ always commutes with direct products. Recall that an $R_2$-module $M$ is called a \emph{cogenerator} of $R_2\Modcat$ if any $R_2$-module can be embedded into a direct product of copies of $M$. Clearly, $I$ is an injective cogenerator of $R_2\Modcat$. In case that $R_2$ is an Artin algebra, there is another injective cogenerator, the usual dual module $D(R_2)$ of the right regular module $R_2$, where $D$ is the usual duality of an Artin algebra.

\begin{Lem}\label{f1}
Let $F:\Dc{R_1}\to\Dc{R_2}$ be a triangle functor. Suppose that $\fd(F)=s>-\infty$ and $\fd(R_2)=t <\infty$. Then we have the following:

$(1)$ $F(\Pf{R_1})\subseteq \mathscr{D}^c_{\geq s-t}(R_2)$.

$(2)$ Let $m\in\mathbb{Z}$. Then, for any $X\in\Pf{R_1}$ and for any $\cpx{Y}\in\D{R_2}$ with $\sup(\cpx{Y})\leq m$, we have $\Hom_{\D{R_2}}(F(X),\cpx{Y}[i])=0$ for all $i>t-s+m$.
\end{Lem}

{\it Proof.} Note that $s=+\infty$ if and only if $F(X)=0$ for any $X\in\Pf{R_1}$. In this case, both $(1)$ and $(2)$ are true. Now, we assume $s<+\infty$. Thus  $s$ is an integer.

(1) Since $F(X)\in\Dc{R_2}$, there exists a complex $\cpx{Q}=(Q^j,d^j)_{j\in\mathbb{Z}}\in\Cb{\pmodcat{R_2}}$ such that $F(X)\simeq \cpx{Q}$ in $\Dc{R_2}$. In particular, $H^i(F(X))\simeq H^i(\cpx{Q})$ for all $i\in\mathbb{Z}$. Since  $\fd(F)=s<\infty$, we have $H^i(F(X))=0$ for all $i<s$. Thus $H^i(\cpx{Q})=0$ for all $i<s$. It follows that $Y:=\Coker(d^{s-1})\in \Pf{R_2}$, and therefore $\cpx{Q}$ is isomorphic in $\D{R_2}$ to the following canonical truncated complex:
$$
0\lra Y\lra Q^{s+1}\lraf{d^{s+1}} Q^{s+2}\lra\cdots\lra 0.
$$
Since $\fd(R_2)=t<\infty$, we have $\pd(_{R_2}Y)\leq t$. So the $R_2$-module $Y$ has a finite projective resolution:
$$
0\lra P^{s-t}\lra\cdots  \lra P^{s-1}\lra P^{s}\lra Y\lra 0
$$
such that $P^j\in\pmodcat{R_2}$ for $s-t\leq j\leq s$. Consequently, $F(X)$ is isomorphic in $\D{R_2}$ to the following complex
$$
\cpx{P}: \quad 0\lra P^{s-t}\lra\cdots  \lra P^{s-1}\lra P^{s}\lra Q^{s+1}\lraf{d^{s+1}} Q^{s+2}\lra\cdots\lra 0.
$$
Clearly, $\cpx{P}\in\Cb{\pmodcat{R_2}}$ and $P^i=0$ for $i<s-t$.
This implies $F(X)\in\mathscr{D}^c_{\geq s-t}(R_2)$. Hence, we have $F(\Pf{R_1})\subseteq \mathscr{D}^c_{\geq s-t}(R_2)$.

$(2)$ Let $X\in\Pf{R_1}$ and $\cpx{Y}\in\D{R_2}$ with $\sup(\cpx{Y})\leq m<\infty$. Then $H^j(\cpx{Y})=0$ for $j>m$, and therefore there exists a complex $\cpx{Z}\in\Cf{R_2}$ with $Z^r=0$ for $r>m$, such that $\cpx{Z}\simeq \cpx{Y}$ in $\D{R_2}$. Moreover, by the proof of $(1)$, there exists another complex $\cpx{P}\in \Cb{\pmodcat{R_2}}$ with $P^i=0$ for all $i<s-t$, such that $\cpx{P}\simeq F(X)$ in $\D{R_2}$. It follows that $$\Hom_{\D{R_2}}(F(X),\cpx{Y}[i])\simeq\Hom_{\D{R_2}}(\cpx{P}, \cpx{Z}[i])\simeq\Hom_{\K{R_2}}(\cpx{P},\cpx{Z}[i])=0$$
for all $i>t-s+m$. This shows $(2)$. $\square$.

\begin{Lem}\label{f2}
Let $F:\Dc{R_1}\to\Dc{R_2}$ be a fully faithful triangle functor. If $\fd(F)=s$ is an integer, then $\fd(R_1)\leq \fd(R_2)-s+\sup(F(R_1))$. In particular, if  $\fd(R_2)<\infty$, then $\fd(R_1)<\infty$.
\end{Lem}

{\it Proof.} If $\fd(R_2)$ is infinity, then the right-hand side of the inequality is infinity and the corollary is true. So we assume that $\fd(R_2)=t<\infty$. Further, we may assume $R_1\ne 0$. Since $F$ is fully faithful, we have $0\neq F(R_1)\in\Dc{R_2}$. This implies that $\sup(F(R_1))<\infty$. Moreover, it is known that, for any $X\in\Pf{R_1}$, if there is a natural number $n$ such that $\Ext^i_{R_1}(X, R_1)=0$ for all $i>n$, then $\pd(_{R_1}X)\leq n$.
So, to show that $\fd(R_1)\leq n:=t-s+\sup(F(R_1))<\infty$, it is enough to prove that
$\Ext^i_{R_1}(X, R_1)=0$ for all $X\in\Pf{R_1}$ and all $i>n$. In fact, since $F$ is fully faithful, we see that $$\Ext^i_{R_1}(X, R_1)\simeq \Hom_{\D{R_1}}(X,R_1[i])\simeq \Hom_{\D{R_2}}(F(X),F(R_1)[i]).$$
Due to Lemma \ref{f1} (2), we have $\Hom_{\D{R_2}}(F(X),F(R_1)[i])=0$ for all $i>n$. Thus $\Ext^i_{R_1}(X, R_1)=0$ for all $X\in\Pf{R_1}$ and all $i>n$. $\square$

\medskip
Summarizing  Lemmas \ref{f0} and \ref{f2} together, we obtain the following useful result, in which $w$ and $cw$ denote the homological width and cowidth of a complex, respectively.

\begin{Koro}\label{f3}
Let $F:\D{R_1}\to \D{R_2}$ be a fully faithful triangle functor such that $F(R_1)\in\Dc{R_2}$.  Then the following statements hold true:

$(1)$ If $F$ has a left adjoint $L:\D{R_2}\to\D{R_1}$ with $L(R_2)\in\Df{R_1}$, then $\fd(R_1)\leq \fd(R_2)+\sup{(L(R_2))}+\sup(F(R_1))$. If moreover $L(R_2)\in\Dc{R_1}$, then $\fd(R_1)\leq \fd(R_2)+ w(L(R_2))$.

$(2)$ If $F$ has a right adjoint $G:\D{R_2}\to\D{R_1}$ and restricts to a functor $\Db{R_1}\to \Db{R_2}$, then $\fd(R_1)\leq \fd(R_2)+cw\big(G(\Hom_\mathbb{Z}(R_2,\mathbb{Q}/\mathbb{Z}))\big)$.
\end{Koro}

{\it Proof.} Clearly, if $\fd(R_2)$ is infinity, then the two statements (1) and (2) are trivially true. So, we assume that $\fd(R_2)=t<\infty$.
We further assume  that $R_i\ne 0$ for $i=1,2.$ By assumption, we have $F(R_1)\in\Dc{R_2}$ , and therefore $F$ restricts to a functor $\Dc{R_1}\to \Dc{R_2}$. Since $F$ is fully faithful and $R_1\neq 0$, we have $F(R_1)\neq 0$. This leads to $\fd(F)\neq +\infty$. Thus $\fd(F)\in \mathbb{Z}\cup\{-\infty\}$.

$(1)$ Since $(L,F)$ is an adjoint pair, we have $H^n(F(R_1)) \simeq \Hom_{\D{R_2}}(R_2, F(R_1)[n]) \simeq \Hom_{\D{R_1}}(L(R_2), R_1[n])$. It follows from $0\neq F(R_1)\in\D{R_2}$ that  $L(R_2)\neq 0$ in $\D{R_1}$. Since $L(R_2)\in \Df{R_1}$, we know that $\sup(L(R_2))$ is an integer. By Lemma \ref{f0}(1), we see that $\inf(F)\geq -\sup(L(R_2))>-\infty$, and therefore $\fd(F)\ge \inf(F)> -\infty$. Combining this with Lemma \ref{f2}, we have
$$\fd(R_1)\leq t-\fd(F)+\sup(F(R_1))\leq t+\sup(L(R_2))+\sup(F(R_1)).$$ This shows the first part of $(1)$. For the second part of $(1)$, we only need to check that $w\big(L(R_2)\big)=\sup{(L(R_2))}+\sup(F(R_1))$.

In fact, it follows from $L(R_2)\in\Dc{R_1}$ that the homological width of $L(R_2)$ is well defined and there exists a complex
$$\cpx{P}: 0\lra P^r \lraf{d^r} P^{r-1}\lra\cdots\lra P^{s-1}\lra P^s\lra 0 $$
in $\Cb{\pmodcat{R_1}}$ with $s=\sup(L(R_2))$ and $s-r= w\big(L(R_2)\big)$ such that $L(R_2)\simeq \cpx{P}$ in $\D{R_1}$ (see Section \ref{sect3.1}). In this case, $d^r$ is not a split injection. Since $(L,F)$ is an adjoint pair, we have $$\Hom_{\D{R_1}}(\cpx{P},R_1[n])\simeq \Hom_{\D{R_1}}(L(R_2),R_1[n])\simeq \Hom_{\D{R_2}}(R_2, F(R_1)[n])\simeq H^n(F(R_1))$$ for all $n\in\mathbb{Z}$. This implies that
$H^n(F(R_1))=0$ for all $n>-r$. Moreover, since  the map $d^r$ is not a split injection, we have $\Hom_{\D{R_1}}(\cpx{P},P^r[-r])\neq 0$. Thus $H^{-r}(F(R_1))\simeq \Hom_{\D{R_1}}(\cpx{P},R_1[-r])\neq 0$.  This shows $\sup(F(R_1))=-r$. It follows that
$w\big(L(R_2)\big)=s-r=\sup{(L(R_2))}+\sup(F(R_1))$.

$(2)$ Under the assumption of $(2)$, we see from Lemma \ref{f0}(2) that
$\inf(F)\geq -\big(m+\id(\Ker(d_{\cpx{I}}^m))\big)$, where $\cpx{I}\in\Cb{\Imodcat{R_1}}$ is defined in Lemma \ref{f0}(2) and $m:=\sup(\cpx{I})$. Thus $\fd(F)\ge \inf(F)>-\infty$ and $$\fd(R_1)\leq t+m+\id(\Ker(d_{\cpx{I}}^m))+\sup(F(R_1))<\infty$$ by Lemma \ref{f2}. Define $\cpx{W}:=G(\Hom_\mathbb{Z}(R_2,\mathbb{Q}/\mathbb{Z}))$.
By the proof of Lemma \ref{f0}(2), we see that $\cpx{W}\simeq \cpx{I}$ in $\D{R_1}$ and that $H^n(\cpx{W})\simeq H^{-n}(F(R_1)$ for all $n\in\mathbb{Z}$.
This implies that $\sup(F(R_1))=-\inf(\cpx{W})=-\inf(\cpx{I})$. Thus $$cw(\cpx{W})=\sup(\cpx{I})-\inf(\cpx{I})+\id(\Ker(d_{\cpx{I}}^m))=m+\sup(F(R_1))+\id(\Ker(d_{\cpx{I}}^m)).$$
So $\fd(R_1)\leq t+cw(\cpx{W})$. $\square$

\medskip
As a consequence of Corollary \ref{f3}, we have the following applicable fact.

\begin{Koro}\label{bimodule}
Let $\cpx{P}\in\C{R_2\otimes_{\mathbb{Z}}R_1\opp}$ such that $_{R_2}\cpx{P}\in\Dc{R_2}$.
Assume that the following conditions hold:

$(1)$ $R_1\simeq\End_{\D{R_2}}(\cpx{P})$ as rings (via multiplication), and $\Hom_{\D{R_2}}(\cpx{P},\cpx{P}[n])=0$ for all $n\neq 0$.

$(2)$ $\cpx{P}_{R_1}$ is isomorphic in $\D{R_1\opp}$ to a bounded complex
$$\cpx{F}: 0\lra F^r \lra F^{r-1}\lra\cdots\lra F^{s-1}\lra F^s\lra 0$$ of flat $R_1\opp$-modules,  where $r,s \in \mathbb{Z}$ and  $r\leq s$.

Then $\fd(R_1)\leq \fd(R_2)+s-r$. In this case, if $\fd(R_2)<\infty$, then $\fd(R_1)<\infty$.
\end{Koro}

{\it Proof.} Let $F:=\cpx{P}\otimesL_{R_1}-:\D{R_1}\to\D{R_2}$.
Then $F(R_1)\simeq {_{R_2}}\cpx{P}\in\Dc{R_2}$ and $F$ has a right adjoint $G:=\rHom_{R_2}(\cpx{P},-):\D{R_2}\to\D{R_1}$. Since $_{R_2}\cpx{P}\in\Dc{R_2}$, the functor $F$ restricts to a functor $F':\Dc{R_1}\to\Dc{R_2}$. Note that the condition $(1)$ implies that $F'$ is fully faithful. Further, since $F$ commutes with direct sums and $\D{R_1}$ is compactly generated by $R_1$, we see that $F$ itself is also fully faithful.

Now, we claim that $F$ restricts to a functor $\Db{R_1}\to\Db{R_2}$. In fact, by Lemma \ref{f0}(2), this is equivalent to saying that  the complex $G\big(\Hom_\mathbb{Z}(R_2,\mathbb{Q}/\mathbb{Z})\big)$ is isomorphic in $\D{R_1}$ to a bounded complex of injective $R_1$-modules.

To  check the latter, we use the functor $(-)^\vee:=\Hom_\mathbb{Z}(-,\mathbb{Q}/\mathbb{Z})$  and apply $G$ to the injective $R_2$-module $R_2^\vee$. Then we have the following isomorphisms in $\D{R_1}$: $$G(R_2^\vee)=\rHom_{R_2}(\cpx{P},R_2^\vee)=\cpx{\Hom}_{R_2}(\cpx{P},R_2^\vee)
\simeq\cpx{\Hom_{\mathbb{Z}}}(R_2\otimes_{R_2}\cpx{P},\mathbb{Q}/\mathbb{Z})\simeq (\cpx{P}){^\vee}.$$
Note that $(-)^\vee:R_1\opp\Modcat\to \,R_1\Modcat$ is an exact functor, which sends flat $R_1\opp$-modules to injective $R_1$-modules. Thus the condition $(2)$ implies that $(\cpx{P}){^\vee}$ is isomorphic in $\D{R_1}$ to the following bounded complex of injective $R$-modules: $$(\cpx{F})^\vee:=0\lra(F^{s})^\vee\lra(F^{s-1})^\vee\lra \cdots\lra  (F^{r-1})^\vee\lra (F^{r})^\vee\lra 0 $$  where $(F^{s})^\vee$ and $(F^{r})^\vee$  are of degrees $-s$ and $-r$, respectively.  Consequently, we have $cw(G({R_2}^\vee))=cw((\cpx{P})^\vee)=cw((\cpx{F})^\vee)\leq s-r$. Now, it follows from Corollary \ref{f3}(2) that
$$\fd(R_1)\leq \fd(R_2)+cw(G({R_2}^\vee))\leq \fd(R_2)+s-r.$$
This completes the proof. $\square$

\medskip
Recall that a ring epimorphism $\lambda:R\to S$ is \emph{homological} if $\Tor^R_i(S,S)=0$ for all $i>0$. This is equivalent to saying that the restriction functor $D(\lambda_*):\D{S}\to \D{R}$
is fully faithful. Note that $D(\lambda_*)$ always has a left adjoint functor $S\otimesL_R-:\D{R}\to\D{S}$. For some new advances on homological ring epimorphisms phrased in terms of recollements of
derived categories, we refer the reader to $\cite{xc3, xc5, xc4}$. Applying Corollary \ref{f3}(1) to homological ring epimorphisms, we have the following simple result.

\begin{Koro}\label{homo-ring}
Let $\lambda:R\to S$ be a homological ring epimorphism such that $_RS\in\Pf{R}$. Then $\fd(S)\leq \fd(R)$. In this case, if $\fd(R)<\infty$, then $\fd(S)<\infty$.
\end{Koro}

{\it Proof.}  If we take $F:=D(\lambda_*)$ and $L:= S\otimesL_R-$ in Corollary \ref{f3}(1), then $\fd(S)\leq \fd(R)+w(L(S))$. Since $w(L(S))=\pd(_SS)=0$, we have $\fd(S)\leq \fd(R).$ $\square$

\smallskip
Let us point out a straightforward proof of Corollary \ref{homo-ring}:

Let $_SX\in\Pf{S}$. Since $\pd(_RS)<\infty$, the Change of Rings Theorem implies that $\pd(_RX)\leq\pd(_SX)+\pd(_RS)<\infty$. Thus $\pd(_RX)\le \fd(R)$. As $\lambda$ is homological, we see that $\Ext_S^i(X,Y)\simeq \Ext_R^i(X,Y)$ for all $Y\in S\Modcat$ and $i\geq 0$.  This implies that $\pd(_SX)\leq \pd(_RX)$. As a result, we have $\pd(_SX)\leq \pd(_RX)\leq\fd(R)$. This shows $\fd(S)\leq \fd(R)$. $\square$

\medskip
The following result extends \cite[Theorem 1.1]{XP} on finitistic dimensions for derived equivalences of left coherent rings to those of arbitrary rings.

\begin{Koro}\label{f3-0}
Suppose that $F:\D{R_1}\to \D{R_2}$ is a triangle equivalence. Then

$$\mid\fd(R_1)-\fd(R_2)\mid\leq w(F(R_1)).$$
\end{Koro}

{\it Proof.} Suppose that $G:\D{R_2}\to\D{R_1}$ is a quasi-inverse of $F$. Then $(G,F)$ and $(F,G)$ are adjoint pairs. Clearly, $G$ is also a triangle equivalence. Since both $F$ and $G$ preserve compact objects, they restrict to triangle
equivalences of perfect derived categories: $F:\Dc{R_1}\lraf{\simeq} \Dc{R_2}$ and $G:\Dc{R_2}\lraf{\simeq}\Dc{R_1}$.
By Corollary \ref{f3}(1), we have $\fd(R_1)\leq \fd(R_2)+w(G(R_2))$ and $\fd(R_2)\leq \fd(R_1)+w(F(R_1))$. Thus, to complete the proof, it is enough to show that $w(G(R_2))=w(F(R_1))$.

In fact, up to isomorphism in derived categories,
we may assume that $F(R_1)\in\Cb{\pmodcat{R_2}}$ and $G(R_2)\in\Cb{\pmodcat{R_1}}$.

Without loss of generality, we suppose that $F(R_1)$ is a complex in $\Cb{\pmodcat{R_2}}$ of the form
$$
0\lra P^{-r} \lraf{d^{-r}} P^{-r+1}\lra\cdots\lra P^{-1}\lra P^0\lra 0
$$
such that $r=w(F(R_1)\geq 0$. This implies that $H^0(F(R_1))\neq 0$ and $d^{-r}$ is not a split injection. Since $(F,G)$ is an adjoint pair, we always have
{\small $$\Hom_{\K{R_2}}(F(R_1),R_2[n])\simeq \Hom_{\D{R_2}}(F(R_1),R_2[n])\simeq \Hom_{\D{R_1}}(R_1,G(R_2)[n])\simeq H^n(G(R_2))$$}for all $n\in\mathbb{Z}$. It follows that $H^i(G(R_2))=0$ for $i<0$ or $i>r$.
Further, we claim that $H^r(G(R_2))\neq 0$, and therefore $\sup(G(R_2))=r$. Actually, since $d^{-r}$ is not a split injection, we have $\Hom_{\K{R_2}}(F(R_1),P^{-r}[r])\neq 0$. Thus $0\neq \Hom_{\K{R_2}}(F(R_1),R_2[r])\simeq H^r(G(R_2))$. So, up to isomorphism in $\K{R_1}$, the complex $G(R_2)$ has the following form
$$
0\lra Q^{s} \lraf{\varphi^{s}} Q^{s+1}\lra Q^{s+2}\lra\cdots
\lra Q^{r-1}\lra Q^r\lra 0\in\Cb{\pmodcat{R_1}}
$$
such that $0\leq r-s=w(G(R_2))$. In particular, this implies that $\varphi^{s}$ is not a split injection. So, to show $w(F(R_1)=w(G(R_2))$, we only need to show $s=0$.

Indeed, since $(G,F)$ is an adjoint pair, we have $$\Hom_{\K{R_1}}(G(R_2),R_1[n])\simeq \Hom_{\D{R_1}}(G(R_2),R_1[n])\simeq \Hom_{\D{R_2}}(R_2,F(R_1)[n])\simeq H^n(F(R_1))$$ for all $n\in\mathbb{Z}.$ On the one hand, if $s<0$, then $\Hom_{\K{R_1}}(G(R_2),R_1[-s])\simeq H^{-s}(F(R_1))=0$, and therefore  $\Hom_{\K{R_1}}(G(R_2),Q^s[-s])=0$. This means that
$\varphi^{s}$ is a split injection, a contradiction. On the other hand, if $s>0$, then
$0=\Hom_{\K{R_1}}(G(R_2),R_1)\simeq H^0(F(R_1))$. This is also a contradiction. Thus $s=0$ and $w(F(R_1))=w(G(R_2))$, as desired.  $\square$

\medskip
The above result describes a relationship for finitistic dimensions of derived equivalent rings. If we weaken derived equivalences into half recollements of perfect derived module categories,
we will obtain the following general result which provides a bound for the finitistic dimension of the middle ring by those of the other two rings.

\begin{Prop}\label{f3-1}
Suppose that there is a half recollement of perfect derived module categories of the rings $R_1, R_2$ and $R_3$
$$
\xymatrix@C=1.2cm{\Dc{R_1}\ar[r]^-{i_*}&\Dc{R_2}\ar[r]^-{j^!}
\ar@/_1.2pc/[l]_-{i^*}&\Dc{R_3}\ar@/_1.2pc/[l]_-{j_!}}\vspace{0.3cm}.
$$
Then $$\fd(R_2)\leq \fd(R_1)+\fd(R_3)+w\big(i_*(R_1)\big)+w\big(j_!(R_3)\big)+1.$$
In particular,
if $\fd(R_1)<\infty$ and $\fd(R_3)<\infty$, then $\fd(R_2)<\infty$.
\end{Prop}

{\bf Proof.} The proof will be done in several steps. We may suppose that $\fd(R_1)<\infty$ and $\fd(R_3)<\infty$. Clearly, if one of $R_1$ and $R_3$ is zero, then Proposition \ref{f3-1} follows from Corollary \ref{f3-0}. From now on, we assume that $R_1\neq 0\neq R_3$.

\smallskip
Step 1. We claim that $j_!j^!(\Pf{R_2})\subseteq\mathscr{D}^{c}_{\geq -u}(R_2)$ where $u:=\fd(R_3)+w\big(j_!(R_3)\big)\geq 0$.

Actually, since $j_!:\Dc{R_3}\to\Dc{R_2}$ is fully faithful, we have
$0\neq j_!(R_3)\in\Dc{R_2}$. This implies $\sup{(j_!(R_3))}<\infty$. As $(j_!,j^!)$ is an adjoint pair, one can follow the proof of Lemma \ref{f0}(1) to show that $-\sup{(j_!(R_3))}\leq \inf(j^!)$. Note that $\inf(j^!)\leq \fd(j^!)$. Thus $-\infty<-\sup{(j_!(R_3))}\leq \fd(j^!)\leq +\infty$.

Define $u_1:=-\sup{(j_!(R_3))}-\fd(R_3)$. Then $u_1\leq \fd(j^!)-\fd(R_3).$ It follows from Lemma \ref{f1}(1) that $$j^!(\Pf{R_2})\subseteq \mathscr{D}^c_{\geq u_1}(R_3).$$
In other words, for any $Y\in\Pf{R_2}$, there exists a complex $\cpx{P}_Y:=(P^n_Y)_{n\in\mathbb{Z}}\in\Cb{\pmodcat {R_3}}$
with $P^n_Y=0$ for $n<u_1$ such that $j^!(Y)\simeq \cpx{P}_Y$ in $\Dc{R_3}$. Clearly, the complex $\cpx{P}_Y$ is of the following form:
$$
0\lra P^{u_1}_Y\lra P^{u_1+1}_Y\lra P^{u_1+2}_Y\lra \cdots\lra P^{s(Y)}_Y \lra 0,
$$
where $s(Y)$ depends on $Y$ and $u_1\leq s(Y)$. Since $j_!(R_3)\in\Dc{R_2}$ by the half recollement, we see that $j_!(R_3)$ is isomorphic in $\Dc{R_2}$ to a bounded complex $\cpx{L}$ of the form
$$
0\lra L^{u_2}\lra L^{u_2+1}\lra L^{u_2+2}\lra \cdots \lra 0
$$
such that $u_2=\sup(j_!(R_3))- w\big(j_!(R_3)\big)$ and that $L^i\in\pmodcat{R_2}$ for all $i\geq u_2$ (see Section \ref{sect3.1}). This implies $j_!(R_3)\in\mathscr{D}^c_{\geq u_2}(R_2)$. Since $\mathscr{D}^c_{\geq u_2}(R_2))$ is closed under direct summands in $\Dc{R_2}$ by Lemma \ref{preparation}(1), we have $j_!(\pmodcat{R_3})\subseteq\mathscr{D}^c_{\geq u_2}(R_2)$.

Note that $u=\fd(R_3)+w\big(j_!(R_3)\big)=\fd(R_3)+\sup{(j_!(R_3))}-u_2=-(u_1+u_2)$. Now, we claim that $j_!j^!(\Pf{R_2})\subseteq\mathscr{D}^c_{\geq -u}(R_2)=\mathscr{D}^c_{\geq (u_1+u_2)}(R_2)$.

Actually, for the complex $\cpx{P}_Y\in\Cb{\pmodcat{R_3}}$, there is a canonical distinguished triangle in $\Dc{R_3}$:
$$
P_Y^{s(Y)}[-s(Y)]\lra \cpx{P}_Y\lra \cpx{P}_Y{^{\leq s(Y)-1}}\lra P_Y^{s(Y)}[1-s(Y)]
$$
where $\cpx{P}_Y{^{\leq s(Y)-1}}$ is truncated from $\cpx{P}_Y$ by replacing $ P^{s(Y)}_Y$ with $0$, that is,
$$
\cpx{P}_Y{^{\leq u_1-1}}:\quad 0\lra P^{u_1}_Y\lra P^{u_1+1}_Y\lra\cdots \lra P^{s(Y)-1}_Y\lra 0\lra 0.
$$
This induces a distinguished triangle in $\Dc{R_2}$:
$$
j_!\big(P_Y^{s(Y)}\big)[-s(Y)]\lra j_!\big(\cpx{P}_Y\big)\lra j_!\big(\cpx{P}_Y{^{\leq s(Y)-1}}\big)\lra j_!\big(P_Y^{s(Y)}\big)[1-s(Y)].
$$
Note that $j_!\big(P_Y^{s(Y)}\big)[-s(Y)]\in \mathscr{D}^c_{\geq s(Y)+u_2}(R_2)\subseteq \mathscr{D}^c_{\geq (u_1+u_2)}(R_2)$ due to $u_1\leq s(Y)$. Since $P_Y^i\in\pmodcat{R_3}$ for all $u_1\leq i\leq s(Y)$, one can apply Lemma \ref{preparation}(2) to show that
$ j_!(\cpx{P}_Y)\in\mathscr{D}^c_{\geq (u_1+u_2)}(R_2)$ by induction on the number of non-zero terms of a complex.  It follows from  $j^!(Y)\simeq\cpx{P}_Y$ that $j_!j^!(Y)\simeq j_!(\cpx{P}_Y)\in\mathscr{D}^c_{\geq(u_1+u_2)}(R_2).$ This implies that $j_!j^!(\Pf{R_2})\subseteq\mathscr{D}^c_{\geq (u_1+u_2)}(R_2)$.

\smallskip
Step 2. We show that $i_*i^*(\Pf{R_2})\subseteq\mathscr{D}^c_{\geq v}(R_2)$, where $v:=\fd(R_1)+w\big(i_*(R_1)\big)+u+1$.

 First of all, we claim that there is an integer $m$ such that $m\leq \fd(i^*)\leq +\infty$. Indeed, the given half recollement yields the following canonical triangle
$$(\dag)\quad
j_!j^!(Y)\lraf{\eta_Y} Y\lraf{\varepsilon_Y} i_*i^*(Y)\lra j_!j^!(Y)[1]$$
in $\D{R_2}$, where $\eta_Y$ and $\varepsilon_Y$ stand for the counit and unit adjunction morphisms, respectively. Since $j_!j^!(Y)\in\mathscr{D}^c_{\geq -u}(R_2)\subseteq\Dc{R_2}$, we can find a complex $\cpx{U}:=(U^n)_{n\in\mathbb{Z}}\in\Cb{\pmodcat{R_2}}$ with $U^n=0$ for all $n<-u\leq 0$ such that $j_!j^!(Y)\simeq\cpx{U}$ in $\D{R_2}$. It follows that $$\Hom_{\D{R_2}}(j_!j^!(Y),Y)\simeq\Hom_{\D{R_2}}(\cpx{U},Y)\simeq\Hom_{\K{R_2}}(\cpx{U},Y).$$
So there exists a chain map $\cpx{f}:\cpx{U}\to Y$ such that its mapping cone $\cpx{V}$ is isomorphic to $i_*i^*(Y)$ in $\D{R_2}$. Clearly, $V^0=U^1\oplus Y$ and $V^j=U^{j+1}$ for any $j\neq 0$. In particular, $V^j=0$ for all $j<-u-1$. Since $i_*:\Dc{R_1}\to \Dc{R_2}$ is fully faithful, we have $$H^n(i^*(Y))\simeq \Hom_{\D{R_1}}(R_1,i^*(Y)[n])\simeq \Hom_{\D{R_2}}(i_*(R_1), i_*i^*(Y)[n])\simeq \Hom_{\D{R_2}}(i_*(R_1), \cpx{V}[n]).$$
By assumption,  $i_*(R_1)\in\Dc{R_2}$ and therefore is isomorphic in $\Dc{R_2}$ to a complex
$$\cpx{Q}:\quad 0\lra Q^{v_2}\lra Q^{v_2+1}\lra\cdots\lra Q^{b} \lra 0
$$
in $\Cb{\pmodcat{R_2}}$, where $b:=\sup(i_*(R_1))$ and $b-v_2=w\big(i_*(R_1)\big)$ (see Section \ref{sect3.1}). Let $m:=-u-1-b$. Then
$$H^n(i^*(Y))\simeq\Hom_{\D{R_2}}(i_*(R_1), \cpx{V}[n])\simeq\Hom_{\D{R_2}}(\cpx{Q}, \cpx{V}[n])\simeq\Hom_{\K{R_2}}(\cpx{Q}, \cpx{V}[n])=0$$
for all $n<m$. This implies that $m\leq \fd(i^*)\leq +\infty$, as claimed.

Let $v_1:=m-\fd(R_1)$. It follows from Lemma \ref{f1}(1) that $i^*(\Pf{R_2})\subseteq \mathscr{D}^c_{\geq v_1}(R_1).$ Now, replacing the pair $(j_!, j^!)$ in the proof of Step $1$ with $(i^*, i_*)$, one can similarly show that $$i_*i^*(\Pf{R_2})\subseteq
\mathscr{D}^c_{\geq v_1+v_2}(R_2).$$
Note that $-(v_1+v_2)=\fd(R_1)+w\big(i_*(R_1)\big)+u+1=v\geq u+1\geq 1$.

Step 3. We show that  $\fd(R_2)\leq v= \fd(R_1)+\fd(R_3)+w\big(i_*(R_1)\big)+w\big(j_!(R_3)\big)+1$.

 Since $j_!j^!(Y)\subseteq\mathscr{D}^c_{\geq -u}(R_2)$ and $i_*i^*(Y)\in\mathscr{D}^c_{\geq -v}(R_2)$ for $Y\in \Pf{R_2}$ with $u< v$, it follows from the triangle $(\dag)$ and Lemma \ref{preparation}(2) that $Y\in\mathscr{D}^c_{\geq -v}(R_2)$. Now, let
$\cpx{P}:=(P^n, d^n)_{n\in\mathbb{Z}}\in\Cb{\pmodcat{R_2}}$ such that $P^n=0$ for all $n<-v$ and that $Y\simeq\cpx{P}$ in $\Dc{R_2}$. Since $Y$ is an $R_2$-module, we see that $H^n(\cpx{P})=0$ for $n\neq 0$ and $H^0(\cpx{P})\simeq Y$. Consequently, $\Ker(d^0)\in\pmodcat{R_2}$ and the following complex
$$
0\lra P^{-v}\lraf{d^{-v}} P^{-v+1}\lraf{d^{-v+1}} \cdots\lra P^{-1}\lraf{d^{-1}} \Ker(d^0)\lra Y\lra 0
$$
is exact. Thus $\pd(_{R_2}Y)\leq v$ and therefore $\fd(R_2)\leq v<\infty$. $\square$

\medskip
Now, with the above preparations, we prove the following strong version of Theorem \ref{main-1} .

\medskip
\begin{Theo}\label{finitistic}
Let $R_1$, $R_2$ and $R_3$ be rings. Suppose that there exists a recollement among
the derived module categories $\D{R_3}$, $\D{R_2}$ and $\D{R_1}$ of $R_3$, $R_2$ and $R_1:$

$$\xymatrix@C=1.2cm{\D{R_1}\ar[r]^-{{i_*}}
&\D{R_2}\ar[r]^-{\;j^!}\ar@/^1.2pc/[l]^-{i^!}\ar_-{i^*}@/_1.2pc/[l]
&\D{R_3} \ar@/^1.2pc/[l]^{j_*\;}\ar@/_1.2pc/[l]_{j_!\;}}$$
Then the following statements hold true:

$(1)$ Suppose that $j_!$ restricts to a functor $\Db{R_3}\to \Db{R_2}$ of bounded derived module categories. Then $\fd(R_3)\leq \fd(R_2)+cw\big(j^!(\Hom_\mathbb{Z}(R_2,\mathbb{Q}/\mathbb{Z}))\big)$.

$(2)$ Suppose that  $i_*(R_1)$ is a compact object in $\D {R_2}$. Then

\;\quad $(a)$ $\fd(R_1)\leq\fd(R_2)+w(i^*(R_2))$.

\;\quad $(b)$ $\fd(R_2)\leq \fd(R_1)+\fd(R_3)+w\big(i_*(R_1)\big)+w\big(j_!(R_3)\big)+1$.
\end{Theo}

{\it Proof.} Note that the triangle functors $j_!$ and $i^*$ in a recollement always take compact objects to compact objects and that $i_*(R_1)$ is compact if and only if $j^!(R_2)$ is compact (for a reference of this fact, one may see, for example, \cite[Lemma 2.2]{xc5}). Thus we have a sequence of functors:
$$
\xymatrix{\Dc{R_1} & \ar[l]_-{i^*}\Dc{R_2} &\ar[l]_-{j_!}\Dc{R_3},}
$$
where the functor $j_!$ is fully faithful.

Applying Corollary \ref{f3}(2) to the adjoint pair $(j_!,j^!)$, we then obtain  $(1)$.

 Suppose $i_*(R_1)\in \Dc{R_2}$. Then $j^!(R_2)\in\Dc{R_3}$ and the given recollement in Theorem \ref{finitistic} induces a half recollment of prefect derived module categories:
 $$
\xymatrix@C=1.2cm{\Dc{R_1}\ar[r]^-{i_*}&\Dc{R_2}\ar[r]^-{j^!}
\ar@/_1.2pc/[l]_-{i^*}&\Dc{R_3}\ar@/_1.2pc/[l]_-{j_!}}\vspace{0.3cm}.$$
Now, the statements $(a)$ and $(b)$ in $(2)$ follow from Corollary \ref{f3}(1) and Proposition  \ref{f3-1}, respectively. This completes the proof of Theorem \ref{finitistic}. $\square$

\medskip
As a consequence of Theorem \ref{finitistic}, we obtain the following corollary which extends the main result \cite[Theorem]{DMX} on finitistic dimensions of Artin algebras to the one of arbitrary rings.

\begin{Koro}\label{stratifying}
Let $R$ be a ring and $e$ an idempotent element of $R$. Suppose that the canonical surjection $R\to R/ReR$ is homological with $_RReR\in\Pf{R}$. Then $\fd(R/ReR)\le \fd(R)\leq \fd(eRe)+\fd(R/ReR)+\pd({_R}R/ReR)+1$.

\end{Koro}

{\it Proof.} Let $J:=ReR$. Since the canonical surjection $R\to R/J$ is homological, there exists a recollement of derived module categories:
$$\xymatrix{\D{R/J}\ar^-{D(\pi_*)}[r]&\D{R}\ar^-{eR\otimes_R^{\mathbb{L}}-}[r]
\ar@/^1.2pc/[l]\ar_-{R/J\otimes_R^{\mathbb
L}-}@/_1.6pc/[l] &\D{eRe}\ar_-{Re\otimes_{eRe}^{\mathbb
L}-}@/_1.6pc/[l]\ar@/^1.2pc/[l].}$$

\medskip
\noindent
Since $_RJ\in\Pf{R}$, we see that $D(\pi_*)(R/J)=R/J\in\Dc{R}$ and that $w(R/J)=\pd(_RR/J)$. Moreover, $Re\otimes_{eRe}^{\mathbb L}eRe=Re$ and $w(_RRe)=0$.
Now, Corollary \ref{stratifying} follows from Theorem \ref{finitistic}(2)(b) and Corollary \ref{homo-ring}. $\square$

\medskip
Since a recollement at $\mathscr{D}^b$-level induces a recollement at $\mathscr{D}$-level, the following result is a straightforward consequence of Theorem \ref{finitistic}, which also generalizes \cite[Theorem 2]{Happel}.

\begin{Koro}\label{db-level}
Let $R_1$, $R_2$ and $R_3$ be rings. Suppose that there exists a recollement among
the derived module categories $\Db{R_3}$, $\Db{R_2}$ and $\Db{R_1}:$

$$\xymatrix@C=1.2cm{\Db{R_1}\ar[r]^-{{i_*}}
&\Db{R_2}\ar[r]^-{\;j^!}\ar@/^1.2pc/[l]^-{i^!}\ar_-{i^*}@/_1.2pc/[l]
&\Db{R_3} \ar@/^1.2pc/[l]^{j_*\;}\ar@/_1.2pc/[l]_{j_!\;}}$$
such that $i_*(R_1)\in \Dc{R_2}$. Then
$$\fd(R_2)<\infty\;\,\mbox{if and only if}\;\,\max{\{\fd(R_1),
\fd(R_3)\}}<\infty.$$
\end{Koro}

\medskip
The existence of a recollement at $\mathscr{D}^b$-level occurs in the following special case (see \cite{NS2}, \cite{koenig}): Let $R$ be a ring and $J=ReR$ be an ideal generated by an idempotent element $e$ in $R$ such that $_RJ$ is projective and finitely generated and that $J_R$ has finite projective dimension. Then there exists a recollement among $\Db{R/J}, \Db{R}$ and $\Db{eRe}$. Remark that, without $\pd(J_R)<\infty$, we may not get a recollement at $\mathscr{D}^b$-level because the left-derived functor
$Re\otimes_{eRe}^{\mathbb L}-:\D{eRe}\to \D{R}$ may not restrict to a functor of bounded derived categories. One can construct a desired counterexample from triangular matrix rings.

Applying Corollary \ref{stratifying} to triangular matrix rings, we re-obtain the following well-known result (for example, see \cite[Corollary 4.21]{FGR}).

\begin{Koro} \label{triangular}
Let $R$ and $S$ be rings, and let $M$ be an $S$-$T$-bimodule. Set $B:=
\left(\begin{array}{lc} S& M\\ 0 & T\end{array}\right).$ Then $\fd(S)\le \fd(B)\leq \fd(S)+\fd(T)+1$.
\end{Koro}

{\it Proof.} Let $e=\left(\begin{array}{lc} 0& 0\\ 0 & 1\end{array}\right)$. Then $BeB=Be$, $eBe\simeq T$, $B/BeB\simeq S=B(1-e)$ and $_BB= {}_BS\oplus Be$. Thus $_BBeB\in\pmodcat{B}$ and the canonical surjection
$B\to S$ is homological. Now, Corollary \ref{triangular} follows from Corollary \ref{stratifying}. $\square$

\medskip
Recall from \cite{xc4} that a morphism $\lambda: Y\ra X$ of objects in an additive category $\mathcal C$ is said to be \emph{covariant}
if the induced map $\Hom_{\mathcal C}(X,\lambda): \Hom_{\mathcal
C}(X,Y)\ra \Hom_{\mathcal C}(X,X)$ is injective, and the induced map $\Hom_{\mathcal C}(Y,\lambda): \Hom_{\mathcal
C}(Y,Y)\ra \Hom_{\mathcal C}(Y,X)$ is a split epimorphism of
$\End_{\mathcal C}(Y)$-modules. Covariant morphisms capture traces of modules, which guarantee the ubiquity of covariant morphisms (see \cite{xc4}).

For covariant morphisms, we have the following result which follows from Corollary \ref{stratifying} and \cite[Lemma 3.2]{xc4}.

\begin{Koro} \label{covariant-dim}
Let $f: Y\ra X$ be a covariant morphism in an additive category $\mathcal{C}$. Then
$$\fd\big(\End_{\mathcal{C},Y}(X)\big)\le\fd\big(\End_{\mathcal{C}}(Y\oplus X)\big)\leq \fd\big(\End_{\mathcal{C}}(Y)\big)+\fd\big(\End_{\mathcal{C},Y}(X)\big)+2,$$
where $\End_{{\mathcal C},Y}(X)$ is the quotient ring of the endomorphism ring
$\End_{\mathcal C}(X)$ of $X$ modulo the ideal generated by all those endomorphisms of $X$ which
factorize through the object $Y$.
\end{Koro}

Consequently, we have the following result which restates Corollary \ref{ars-dim}.

\begin{Koro} \label{ars-idem}
($1$) Let $I$ be an idempotent ideal in a ring $R$. Then
$$\fd(R/I)\le\fd\big(\End_R(R\oplus I)\big)\leq \fd\big(\End_R(_RI)\big)+\fd(R/I)+2.$$
In particular, if $_RI$ is projective and finitely generated, then $$\fd(R/I)\le\fd(R)\leq \fd\big(\End_R(_RI)\big)+\fd(R/I)+2.$$

$(2)$ Let $0\ra Z\ra Y\lraf{f} X\ra 0$ be an almost split sequence in $R\modcat$ with $R$ an Artin algebra such that $\Hom_R(Y,Z)=0$ (see \cite{ars} for definition). Then $$\fd\big(\End_R(Y\oplus X)\big)\leq \fd\big(\End_R(Y)\big)+2.$$

\end{Koro}

{\it Proof.} (1) Since the inclusion  $I\hookrightarrow R$ is a covariant homomorphism in $R\Modcat$ and $\End_{R,I}(R)\simeq R/I$, we know that the first statement in (1) follows from Corollary \ref{covariant-dim} immediately. The last statement is a consequence of the fact that $R$ is Morita equivalent to $\End_R(R\oplus I)$.

(2) Under the assumption, we know that $f$ is a covariant map in $R\modcat$, the category of finitely generated $R$-modules.
So, by Corollary \ref{covariant-dim}, it is sufficient to show that $\fd(\End_{R,Y}(X))=0$. In fact, since $\End_R(X)$ is a local algebra and since the ideal of $\End_R(X)$ generated by all homomorphisms which factorize through $Y$ belong to the radical of $\End_R(X)$, the algebra $\End_{R,Y}(X)$ is local. Note that a local Artin algebra has finitistic dimension $0$. Therefore $\fd(\End_{R,Y}(X))=0$. Now, (2) follows from Corollary \ref{covariant-dim}. $\square$

Note that an alternative proof of Corollary \ref{ars-dim}(2) can be given by \cite[Theorem 1.1]{hx2} together with Corollary \ref{triangular} and \cite{XP}.

\medskip
In the following we point out that the methods developed in this paper for little finitistic dimensions also work for big finitistic and global dimensions. Recall that, for an arbitrary ring $R$, we denote by $\Fd(R)$ and $\gd(R)$ the \emph{big finitistic} and \emph{global dimensions} of $R$, respectively. By definition, $\gd(R)$ (respectively, $\Fd(R)$) is the supremum of projective dimensions of all left $R$-modules (respectively, which have finite projective dimension). Clearly, $\fd(R)\leq \Fd(R)\leq \gd(R)$; and if $\gd(R)<\infty$, then $\Fd(R)=\gd(R)$. However, the equality $\fd(R)=\Fd(R)$ does not have to hold in general (see \cite{ZH}).

As in Theorem \ref{finitistic}, we have the following result on big finitistic dimensions of rings, in which the condition (2) is  weaker than the one in Theorem \ref{finitistic}(2).

\begin{Theo}\label{big finitistic}
Let $R_1$, $R_2$ and $R_3$ be rings. Suppose that there exists a recollement among
the derived module categories $\D{R_3}$, $\D{R_2}$ and $\D{R_1}$ of $R_3$, $R_2$ and $R_1:$

$$\xymatrix@C=1.2cm{\D{R_1}\ar[r]^-{{i_*}}
&\D{R_2}\ar[r]^-{\;j^!}\ar@/^1.2pc/[l]^-{i^!}\ar_-{i^*}@/_1.2pc/[l]
&\D{R_3} \ar@/^1.2pc/[l]^{j_*\;}\ar@/_1.2pc/[l]_{j_!\;}}$$
Then the following statements hold true:

$(1)$ Suppose that $j_!$ restricts to a functor $\Db{R_3}\to \Db{R_2}$ of bounded derived module categories. Then $\Fd(R_3)\leq \Fd(R_2)+cw\big(j^!(\Hom_\mathbb{Z}(R_2,\mathbb{Q}/\mathbb{Z}))\big)$.

$(2)$ Suppose that  $i_*(R_1)$ is isomorphic in $\D{R_2}$ to a bounded complex of (not necessarily finitely generated) projective $R_2$-modules. Then we have the following:

\;\quad $(a)$ $\Fd(R_1)\leq\Fd(R_2)+w(i^*(R_2))$.

\;\quad $(b)$ $\Fd(R_2)\leq \Fd(R_1)+\Fd(R_3)+w\big(i_*(R_1)\big)+w\big(j_!(R_3)\big)+1$.
\end{Theo}

{\it Sketch of the proof.} Let us consider the full subcategory $\mathscr{X}(R)$ of $\D{R}$ consisting of all those complexes which are isomorphic in $\D{R}$ to bounded complexes of projective $R$-modules. It is known that $\mathscr{X}(R)$ contains $\Dc{R}$ and that the localization functor $\K{R}\to\D{R}$ induces a triangle equivalence $\Kb{\Pmodcat{R}}\lraf{\simeq}\mathscr{X}(R)$.

Similarly, one can define big finitistic dimensions of functors, and establish several parallel results for $\Fd(R)$, such as Lemma \ref{f2}, Corollary \ref{f3} and Proposition \ref{f3-1}. In the present situation, we shall replace $\Dc{R}$ with $\mathscr{X}(R)$, and consider big finitistic dimensions of triangle functors which commute with direct sums. Further, to show Theorem \ref{big finitistic}, we observe the following facts for a given recollement:

$(i)$ The functors $j_!$, $j^!$, $i^*$ and $i_*$ commute with direct sums.

$(ii)$ The functors $j_!$ and $i^*$ preserve compact objects and restrict to triangle functors $$\mathscr{X}(R_3)\lraf{j_!} \mathscr{X}(R_2)\;\;\mbox{and}\;\;\mathscr{X}(R_2)\lraf{i^*} \mathscr{X}(R_1).$$

$(iii)$ If $i_*(R_1)\in\mathscr{X}(R_2)$, then $i_*$ and $j^!$ restrict to triangle functors $$\mathscr{X}(R_1)\lraf{i_*} \mathscr{X}(R_2)\;\;\mbox{and}\;\;\mathscr{X}(R_2)\lraf{j^!} \mathscr{X}(R_3).$$
Now, one can use the methods in the proof of Theorem \ref{finitistic} to show Theorem
\ref{big finitistic}. Here, we omit the details. $\square$

\medskip
Concerning global dimensions, we can describe explicitly upper bounds for the global dimension of a ring in terms of the ones of the other two rings involved in a recollement. These upper bounds imply the finiteness of global dimensions mentioned in \cite[Proposition 2.14]{AKLY}.

\begin{Theo}\label{gldim}
Let $R_1$, $R_2$ and $R_3$ be rings. Suppose that there exists a recollement among
the derived module categories $\D{R_3}$, $\D{R_2}$ and $\D{R_1}$ of $R_3$, $R_2$ and $R_1:$

$$\xymatrix@C=1.2cm{\D{R_1}\ar[r]^-{{i_*}}
&\D{R_2}\ar[r]^-{\;j^!}\ar@/^1.2pc/[l]^-{i^!}\ar_-{i^*}@/_1.2pc/[l]
&\D{R_3} \ar@/^1.2pc/[l]^{j_*\;}\ar@/_1.2pc/[l]_{j_!\;}}$$
Then we have the following:

$(1)$ If $\gd(R_2)<\infty$, then $\gd(R_3)\leq \gd(R_2)+cw\big(j^!(\Hom_\mathbb{Z}(R_2,\mathbb{Q}/\mathbb{Z}))\big)$ and $\gd(R_1)\leq\gd(R_2)+w(i^*(R_2))$.

$(2)$  If $\gd(R_1)<\infty$ and $\gd(R_3)<\infty$, then $\gd(R_2)\leq \gd(R_1)+\gd(R_3)+w\big(i_*(R_1)\big)+w\big(j_!(R_3)\big)+1$.
\end{Theo}

{\it Sketch of the proof.} From \cite[Proposition 2.14]{AKLY} and its proof, we observe the following two facts:

$(i)$ If $\gd(R_2)<\infty$ or $\gd(R_3)<\infty$, then $i_*(R_1)$ is isomorphic in $\D{R_2}$ to a bounded complex of projective $R_2$-modules.

$(ii)$ $\gd(R_2)<\infty$ if and only if both $\gd(R_1)<\infty$ and $\gd(R_3)<\infty$. In this case, the recollement among unbounded derived categories can restrict to a recollement of bounded derived categories.

Moreover, for a ring $R$, if $\gd(R)<\infty$, then $\gd(R)=\Fd(R)$. Now, Theorem \ref{gldim} becomes a consequence of Theorem \ref{big finitistic}. $\square$

\subsection{Proofs of Theorem \ref{hom-dim} and Corollary \ref{ring extension}}
Now let us turn to proofs of our results on exact contexts that arise from different situations.

\medskip
{\bf Proof of Theorem \ref{hom-dim}.}

Given an exact context $(\lambda,\mu, M, m),$ we have defined its noncommutative tensor
product $T\boxtimes_RS$ and
the following two ring homomorphisms
$$
\rho: S\to T\boxtimes_RS,\; s\mapsto 1\otimes s\;\;\mbox{for}\;\;s\in S,\;\;\mbox{and}\;\;
\phi: T\to T\boxtimes_RS,\; t\mapsto t\otimes 1\;\;\mbox{for}\;\; t\in
T.
$$
Note that $T\boxtimes_RS$ has $T\otimes_RS$ as its abelian groups, while its multiplication is different from the usual tensor product (see \cite{xc3} for details). Let $B:=\left(\begin{array}{lc} S& M\\ 0 & T\end{array}\right)$, $C:=M_2(T\boxtimes_RS)$ and
$$ \theta:=\left(\begin{array}{cc} \rho & \beta\\
0 & \phi\end{array}\right): \; B \lra C, $$ where $\beta: M\ra T\otimes_RS$ is the unique $R$-$R$-bimodule homomorphism such that $\phi=(m\cdot)\beta$ and $\rho=(\cdot m)\beta$.

Let $$
\varphi: \left(\begin{array}{l} S \\
0 \end{array}\right)\lra \left(\begin{array}{l} M \\
T\end{array}\right),\;\left(\begin{array}{l} s \\
0 \end{array}\right)\mapsto \left(\begin{array}{c} sm\\
0 \end{array}\right)\;\;\mbox{for}\;\; s\in S.$$
Then $\varphi$ is a homomorphism of $B$-$R$-bimodules. Denote by $\cpx{P}$ the mapping cone
of $\varphi$. Then $\cpx{P}\in\Cb{B\otimes_{\mathbb{Z}}R\opp}$ and $_B\cpx{P}\in\Cb{\pmodcat{B}}$.
In particular, $\cpx{P}\in\Dc{B}$.

By \cite[Theorem 1.1]{xc3}, if $\Tor^R_i(T,S)=0$ for all $i\geq 1$, then  there is a recollement of derived categories:
$$
\xymatrix@C=1.2cm{\D{C}\ar[r]^-{D(\theta_*)}
&\D{B}\ar[r]^-{j^!}\ar@/^1.2pc/[l]\ar_-{C\otimesL_B-}@/_1.2pc/[l]
&\D{R} \ar@/^1.2pc/[l]\ar@/_1.2pc/[l]_{j_!\;}}$$

\medskip
\noindent where $j_!:={_B}\cpx{P}\otimesL_R-$, $j^!:=\cpx{\Hom}_B(\cpx{P},-)$ and $D(\theta_*)$ is the restriction functor induced from the ring homomorphism $\theta:B\to C$. First of all, we have two easy observations:

(i) Since $C:=M_2(T\boxtimes_RS)$ is Morita equivalent to $T\boxtimes_RS$, we have $\fd(C)=\fd(T\boxtimes_RS)$.

(ii) Since $B$ is a triangular matrix ring with the rings $S$ and $T$ in the diagonal, it follows from Corollary \ref{triangular} that $\fd(B)\leq \fd(S)+\fd(T)+1$.

\smallskip
We first apply Corollary \ref{bimodule} to show Theorem \ref{hom-dim}(1). In fact, by
\cite[Lemma 5.4]{xc3}, we see that $R\simeq \End_{\D{B}}(\cpx{P})$ as rings (via multiplication) and that $\Hom_{\D B}\big(\cpx{P},\cpx{P}[n]\big)=0$ for any $n\neq 0$. It remains to show that
$\cpx{P}_R$ is isomorphic in $\D{R\opp}$ to a bounded complex of flat $R\opp$-modules.

Since the sequence $ 0\ra R\lraf{(\lambda,\,\mu)}S\oplus T \lraf{\left({\cdot\,m\,\atop{-m\,\cdot}}\right)}M\ra 0 $
is exact, we have $\Cone(\cdot m)\simeq \Cone(\mu)$ in $\D{R\opp}$. This implies that $\cpx{P}_R\simeq T\oplus \Cone(\cdot m)\simeq T\oplus \Cone(\mu)$ in $\D{R\opp}$, where $\Cone(\mu)$ is the complex $0\ra R\lraf{\mu} T\ra 0$ with $T$ in degree $0$.
If $\fld(T_R)=\infty$, then Theorem \ref{hom-dim}(1) is trivially true. So we may suppose $\fld(T_R)<\infty$. Let $t:=\max\{1,\fld(T_R)\}$. Then $\cpx{P}$ is isomorphic in $\D{R\opp}$ to a bounded complex $$\cpx{F}:=0\lra F^{-t}\lra F^{-t+1}\lra \cdots\lra F^{-1}\lra F^{0}\lra 0 $$ such that $F^i$ are flat $R\opp$-modules for $-t\leq i\leq 0$. It follows from Corollary \ref{bimodule} that  $\fd(R)\leq \fd(B)+t\leq \fd(S)+\fd(T)+t+1.$ This shows Theorem \ref{hom-dim}(1).

Next, we shall apply Theorem \ref{finitistic} to the above recollement and give a proof of Theorem \ref{hom-dim}(2).

By the proof of \cite[Theorem 1.3(2)]{xc5}, we see that $D(\theta_*)(C)={_B}C\in\Pf{B}$ if and only if $_RS\in\Pf{R}$. Suppose $_RS\in\Pf{R}$. It follows from \cite[Corollary 5.8(1)]{xc3} that $$\pd(_BC)\leq \max\{2,\pd(_RS)+1\}.$$

Since $C\otimesL_BB\simeq C$ in $\D{C}$, we see from Theorem \ref{finitistic}(2)(a) that $\fd(C)\leq \fd(B)$. Note that $\fd(T\boxtimes_RS)=\fd(C)$ and $\fd(B)\leq \fd(S)+\fd(T)+1$.
Thus $(a)$ holds.

Since $D(\theta_*)(C)={_B}C$ and $j_!(R)\simeq {_B}\cpx{P}$ in $\D{B}$, we know that $w(\cpx{P})=1$ and $$w\big(D(\theta_*)(C)\big)=w(_BC)=\pd(_BC)\leq\max\{2,\pd(_RS)+1\}.$$   Now, it follows from Theorem \ref{finitistic}(2)(b) that $$\fd(B)\leq \fd(R)+\fd(T\boxtimes_RS)+\max\{2,\,\pd(_RS)+1\}+1+1.$$ Clearly, $\fd(S)\leq \fd(B)$. Thus $(b)$ holds. $\square$

\medskip
Let us point out the following fact related to Theorem \ref{hom-dim}(2): Suppose that $(\lambda, \mu, M,m)$ is an exact content with $\Tor_i^R(T,S)=0$ for all $i\geq 1$. If $\lambda:R\to S$ is a homological ring epimorphism such that $_RS\in \Pf{R}$, then $\fd(S)\leq \fd(R)$ and $\fd(T\boxtimes_RS)\leq \fd(T)$.

In fact, in this case, the Tor-vanishing condition, that is, $\Tor^R_i(T,S)=0$ for all $i>0$, is equivalent to that $\phi:T\to T\boxtimes_RS$ is a homological ring epimorphism (see \cite[Theorem 1.1(1)]{xc3} for details). Moreover, we have $T\boxtimes_RS\simeq T\otimes_RS$ as $T$-$S$-bimodules. It follows that if $_RS\in\Pf{R}$, then $_TT\boxtimes_RS\in\Pf{T}$ by the Tor-vanishing condition. Therefore the above-mentioned fact is a consequence of Corollary \ref{homo-ring}.

\medskip
{\bf Proof of Corollary \ref{ring extension}.}

Let $\tau: S\subseteq R$ be the inclusion of from $S$ into $R$, and let $\pi:R\to R/S$ be the canonical surjection. We define $$\sigma:\,S \lra R'=\End_S(R/S), \quad s\mapsto (r\mapsto rs)\;\,\mbox{for}\;s\in S\;\mbox{and}\; r\in R/S$$ to be the right multiplication map. Then the quadruple $\big(\tau, \sigma, \Hom_S(R,R/S), \pi\big)$ determined by the extension is an exact context (see the examples in \cite[Section 3]{xc3}) and its noncommutative tensor product $R'\boxtimes_S R$ is defined. If $_SR$ is flat, then $\Tor^S_i(R', R)=0$ for all $i\geq 1$. Particularly, under the assumption on $_SR$ in Corollary \ref{ring extension}(2), the quadruple fulfills the Tor-vanishing condition in Theorem \ref{hom-dim}(2).

Now, we apply Theorem \ref{hom-dim} to the exact context $(\tau, \sigma, \Hom_S(R,R/S), \pi)$, and see that the statements $(a)$ and $(b)$ in Corollary \ref{ring extension} follow from the statements $(a)$ and $(b)$ in Theorem \ref{hom-dim}, respectively. To show  Corollary \ref{ring extension}(1), we shall apply Theorem \ref{hom-dim}(1). For this aim,  we shall prove $$\fld(R'_S)\leq \max{\{\fld\big(\Hom_S(R,R/S)_S\big), \fld\big((R/S)_S\big)\}}.$$
However, this can be concluded from the following exact sequence of right $R'$-modules (also right $S$-modules):
$$
0\lra R'\lra \Hom_S(R,R/S)\lra \Hom_S(S,R/S)\lra 0.
$$
which is obtained by applying $\Hom_S(-,R/S)$ to the exact sequence $0\to S\to R\to R/S\to 0 $.
Now, the statement $(1)$ follows from Theorem \ref{hom-dim}(1). $\square$

\subsection{Proofs of Corollaries \ref{mod1b} and \ref{mod1a}}

In the following, we shall show that under the assumptions in Corollaries \ref{mod1b} and \ref{mod1a}, we can get exact pairs, a special class of exact contents, which satisfy the Tor-vanishing condition in Theorem \ref{hom-dim}, and then apply Theorem \ref{hom-dim} to each case. Here, noncommutative tensor products will be replaced by coproducts, and the latter will be interpreted further as some usual constructions of rings.

Let $\lambda: R\ra S$ and $\mu: R\ra T$ be ring homomorphisms, and let $M$ be an $S$-$T$-bimodule with $m\in M$. Recall that an exact context $(\lambda, \mu, M, m)$ is called an \emph{exact pair} if $M=S\otimes_RT$ and $m=1\otimes 1$. In this case, we simply say that $(\lambda,\mu)$ is an exact pair. By \cite[Corollary 4.3]{xc3}, if the map $\lambda$ in the exact context
is a ring epimorphism, then the pair $(\lambda,\mu)$ is exact.
Moreover, by \cite[Remark 5.2]{xc3}, for an exact pair $(\lambda,\mu)$, we have $T\boxtimes_RS\simeq S\sqcup_RT$, the coproduct of the $R$-rings of $S$ and $T$.

\medskip
{\bf Proof of Corollary \ref{mod1b}}.

We define $T:=R\ltimes M$, $\mu: R\to T$ to be
the inclusion from $R$ into $T$, and $\widetilde{\lambda}: R\ltimes M\ra S\ltimes M$ to be the canonical map induced from $\lambda$. By Lemma \ref{tri-ext}, the ring
$S\ltimes M$, together with the inclusion $\rho: S\ra S\ltimes M$
and $\widetilde{\lambda}: T\to S\ltimes M$, is the coproduct of $S$
and $T$ over $R$.

Now, we show that $(\lambda, \mu)$ is an exact pair.
Actually, the split exact sequence $0\ra
R\lraf{\mu} T\ra M\ra 0$  of $R$-$R$-bimodules implies that $_RT_R\simeq
R\oplus M$ as $R$-$R$-bimodules. Since $\lambda$ is a ring
epimorphism and $M$ is an $S$-$S$-bimodule, the map
$$S\otimes _R T\lra  S\ltimes M, \; s\otimes (r, m)\mapsto (sr,
sm)$$ for $s\in S$ and $m\in M$, is an isomorphism of
$S$-$T$-bimodules. Under this isomorphism, we can identify the map
$\mu\,'=id_S\otimes \mu: S\to  S\otimes T$ with the inclusion $\rho: S\ra S\ltimes M$, and the map $\lambda'=\lambda\otimes id_T: T\ra S\otimes_RT$ with $\widetilde{\lambda}$.  Note that $0\ra S \lraf{\rho} S\ltimes M \ra M\ra
0$ is also a split exact sequence of $S$-$S$-bimodules. It follows
that $\Coker(\mu)\simeq\Coker(\rho)\simeq M$ as $R$-$R$-bimodules,
and therefore the sequence of $R$-$R$-bimodules:
$$
0\lra R\lraf{(\lambda,\,\mu)}S\oplus T
\lraf{\left(\rho\atop{-\widetilde{\lambda}}\right)} S\ltimes M\lra 0
$$is exact. This means that the pair $(\lambda, \mu)$ is exact.

Consequently, we know that $T\boxtimes_RS\simeq S\sqcup_RT\simeq S\ltimes M$ as rings.
Note that $\Tor^R_i(T, S)\simeq \Tor^R_i(R\oplus M, S)\simeq
\Tor^R_i(M, S)=0$ for all $i\geq 1$. Thus Corollary \ref{mod1b}(a) follows immediately from Theorem \ref{hom-dim}(2)(a).

Now we turn to the proof of Corollary \ref{mod1b}(b).

Note that, if we apply Theorem \ref{hom-dim}(2)(b) to the exact pair $(\lambda,\mu)$, then we only get $\fd(S)\leq \fd(R)+\fd(S\ltimes M)+\max\{1,\pd(_RS)\}+3$. So, to obtain the better upper bound given in Corollary \ref{mod1b}(b), we need the following statement:

\smallskip
$(\ast)\;\;$
Let $f:\Lambda\to\Gamma$ and $g:\Gamma\to\Lambda$ be ring homomorphisms such that $fg={\rm Id}_{\Lambda}$. If $\fd({_\Gamma}\Gamma\otimesL_\Lambda-)=s< +\infty$, then
$\fd(\Lambda)\leq \fd(\Gamma)-s$.

\smallskip
To show $(*)$, we set $F:={_\Gamma}\Gamma\otimesL_\Lambda-:\D{\Lambda}\to \D{\Gamma}$. If $\fd(F)=-\infty$, then $(*)$ is automatically true. So, we suppose  that $\fd(F)=s$ is an integer and $\Lambda\neq 0$. Since $F(\Lambda)\simeq \Gamma\neq 0$, we have $s\leq 0$. Let $X\in\Pf{\Lambda}$. Then there exists a finite projective resolution  $0\ra P_n\ra \cdots \ra P_1\ra P_0\ra X\ra 0$ of ${_\Lambda}X$ with all $P_i$ in $\pmodcat{\Lambda}$. Now we define $Y:=\Omega^{-s}_\Lambda(X)$, the $(-s)$-th syzygy module of ${_\Lambda}X$. Thus $Y\in\Pf{\Lambda}$. Since $\fd(F)=s$, we see that  $\Tor_{j}^\Lambda(\Gamma,Y)= \Tor_{j-s}^\Lambda(\Gamma,X)\simeq H^{s-j}(F(X))=0$ for all $j>0$. It follows that $\Gamma\otimes_\Lambda Y\in\Pf{\Gamma}$ and $\Gamma\otimes_\Lambda\Omega^i_\Lambda(Y)=\Omega^i_\Gamma(\Gamma\otimes_\Lambda Y)\oplus Q_i$ for all $i\geq 0$, where all $Q_i$ are finitely generated projective $\Gamma$-modules. Further, we may suppose that $\fd(\Gamma)=t<\infty$. Then $\pd(_\Gamma \Gamma\otimes_\Lambda Y)\leq t$, and therefore $\Gamma\otimes_\Lambda\Omega^t_\Lambda(Y) =\Omega^t_\Gamma(\Gamma\otimes_\Lambda Y)\oplus Q_t\in\pmodcat{\Gamma}$.
Since $fg={\rm Id}_{\Lambda}$, we have $\Omega^t_\Lambda(Y)\simeq \Lambda\otimes_\Gamma(\Gamma\otimes_\Lambda\Omega^t_\Lambda(Y))\in\pmodcat{\Lambda}$.
Consequently,  $$\pd(_\Lambda X)\leq \pd(_\Lambda Y)-s\leq \pd(_\Lambda\Omega^t_\Lambda(Y))+t-s\leq t-s.$$ Thus $\fd(\Lambda)\leq \fd(\Gamma)-s$. This finishes the proof of $(\ast)$.

Now, we take $\Lambda:=S$ and $\Gamma:=S\ltimes M$.  Let $f:S\to S\ltimes M$ and $g:S\ltimes M\to S$ be the canonical injection and surjection, respectively. Clearly, we have $fg=\rm{Id}_S$. We assume $S\neq 0$. Then $_{\Gamma}\Gamma\otimesL_{\Lambda}\Lambda=\Gamma \ne 0$, and $\fd(_{\Gamma}\Gamma\otimesL_{\Lambda}-)\le 0$.
Suppose $\fd(R)=m<\infty$. Due to $(\ast)$, in order to show Corollary \ref{mod1b}(b), we only need to prove that $\fd({_\Gamma}\Gamma\otimesL_S-)\geq -m$. This is equivalent to saying that $\Tor^S_n(\Gamma,X)\simeq\Tor^S_n(M,X)=0$ for all $X\in\Pf{S}$ and for all $n>m$.

To check the latter, we first prove that $\Tor^S_j(M,N)\simeq \Tor^R_j(M,N)$ for any $S$-module $N$ and for all $j\geq 1$. Indeed, let $\cpx{P}$ be a deleted projective resolution of the $R\opp$-module $M$. Since   $\Tor^R_i(M,S)=0$ for all $i\geq 1$, we see that $\cpx{P}\otimes_RS$ is a deleted projective resolution of the $S\opp$-module $M\otimes_RS$. Note that $M\otimes_RS\simeq M$ as $S\opp$-modules since $\lambda:R\to S$ is a ring epimorphism and $M$ is an $S\opp$-module. It follows that $\cpx{P}\otimes_RS$ is a deleted projective resolution of the $S\opp$-module $M$. Since $(\cpx{P}\otimes_RS)\otimes_SN\simeq \cpx{P}\otimes_RN$ as complexes, we have $\Tor^S_j(M,N)\simeq \Tor^R_j(M,N)$ for all $j\geq 1$.

Let $_SX\in\Pf{S}$. Since $\pd(_RS)<\infty$, the Change of Rings Theorem implies that $\pd(_RX) \leq \pd(_SX) + \pd(_RS) < \infty$. Hence $\pd(_RX)\le m=\fd(R)$ and $\Tor^S_n(M,X)\simeq \Tor^R_n(M,X)=0$ if $n>m$. This implies that $\fd({_\Gamma}\Gamma\otimesL_S-)\geq -m$.
Now, by the result $(\ast)$, we obtain $\fd(S)\leq \fd(\Gamma)+m=\fd(S\ltimes M)+\fd(R).$
This completes the proof of Corollary \ref{mod1b}(b). $\square$

\medskip
We remark that the statement $(\ast)$ also implies that for any trivial extension of $R$ by an $R$-$R$-bimodule $M$, we always have $\fd(R)\leq \fd(R\ltimes M) +\fld(M_R)$.

\medskip
{\bf Proof of Corollary \ref{mod1a}}:

Let $\lambda:R\to S:=R/I_1$ and $\mu:R\to T:=R/I_2$ be the canonical surjective ring homomorphisms. Since $I_1\cap I_2=0$, we see that $(\lambda,\mu,R/(I_1+I_2),1)$ is an exact context, where $1$ is the identity of the ring $R/(I_1+I_2)$. Even more, since $R$ is a pullback of the surjective maps $R\ra R/I_i$ over $R/(I_1+I_2)$, the pair $(\lambda,\mu)$ is
exact (for example, see \cite[Section 3]{xc3}). So $T\boxtimes_RS\simeq S\sqcup_RT$ as rings.
Note that $S\sqcup_RT=(R/I_1)\sqcup_R(R/I_2)=R/(I_1+I_2)$ by Lemma \ref{epi00'}(2). Thus $T\boxtimes_RS\simeq R/(I_1+I_2)$ as rings.

Now, we apply Theorem \ref{hom-dim} to show Corollary \ref{mod1a}. Clearly, it remains to check that if $\Tor_i^R(I_2,I_1)=0$ for $i\geq 0$, then $\Tor^R_i(R/I_2,R/I_1)=0$ for all $i>0$. In fact, for $i>2$, we have $\Tor^R_i(R/I_2,R/I_1)\simeq\Tor^R_{i-2}(I_2,I_1)=0$ by assumption. Note that $\Tor^R_1(R/I_2,R/I_1)\simeq (I_2\cap I_1)/(I_2I_1)=0$ and $\Tor^R_2(R/I_2,R/I_1)\simeq \Tor^R_1(I_2,R/I_1)=\Ker(f)$ where $f:I_2\otimes_RI_1\to I_2I_1$ is the multiplication map.
Since $I_2\otimes_RI_1=0$, we have $\Tor^R_2(R/I_2,R/I_1)\simeq \Ker(f)=0$. Thus $\Tor^R_i(R/I_2,R/I_1)=0$ for all $i>0$.  $\square$

\medskip
Finally, we apply our results to exact contexts related to homological ring epimorphisms.
First of all, we establish a method to
construct new homological ring epimorphisms from given ones.

\begin{Lem}\label{add}
Let $\lambda: R\to S$ be a homological ring epimorphism. Suppose
that $I$ is an ideal of $R$ such that the image $J'$ of $I$ under
$\lambda$ is a left ideal in $S$ and that the restriction of
$\lambda$ to $I$ is injective. Let  $J$ be the ideal of $S$
generated by $J'$.  Then the following statements are
equivalent:

$(1)$ The homomorphism $\widetilde{\lambda}: R/I\ra S/J$ induced from $\lambda$ is
homological.

$(2)$ $\Tor^{R/I}_i\big(J/J',\,S/J\big)=0$ for all $i\geq 1$.

$(3)$ The multiplication map $I\otimes_RS\to J$ is an isomorphism
and $\Tor^{R}_j(I,\,S)=0$ for all $j\geq 1$.

$(4)$ $\Tor^R_j(R/I,S)=0$ for all $j\ge 1$.

\medskip
Let $B:=\left(\begin{array}{lc} S & S/J'\\
0 & R/I\end{array}\right)$. If one of the above statements holds true, then there is a
recollement of derived module categories:
$$
\xymatrix@C=1.2cm{\D{S/J}\ar[r]&\D{B}\ar[r]
\ar@/^1.2pc/[l]\ar@/_1.2pc/[l]
&\D{R}\ar@/^1.2pc/[l]\ar@/_1.2pc/[l]}.
$$\vspace{0.1cm}
\end{Lem}

{\it Proof.} We take $T:=R/I$ and choose $\mu:R\to T$ to be the canonical surjective homomorphism of rings. Since $J'$ is a left ideal of $S$, we have $S\otimes_RT=S\otimes_R(R/I)\simeq S/(S\cdot I)=S/J'$. On the one hand, the pair $(\lambda,\mu)$ is exact if and only if $\lambda|_I:I\to J'$ is an isomorphism. On the other hand, by Lemma \ref{epi00'}(2), we see that $S\sqcup_{R}T=S\sqcup (R/I) =S/J$ with $J=J'S$, and that the ring homomorphism  $\phi:T\ra
S\sqcup_R T$ in \cite[Theorem 1.1]{xc3} can be chosen as the canonical
map $\widetilde{\lambda}: R/I\to S/J$ induced from $\lambda$.
Thus $(1)$ and $(4)$ are equivalent by \cite[Theorem 1.1(1)]{xc3}. Moreover, the recollement
follows from \cite[Theorem 1.1(2)]{xc3}.

In the following, we shall show that $(3)$ and $(4)$ are equivalent.

Applying the tensor functor $-\otimes_RS$ to the exact sequence $0\ra I\ra R\ra R/I\ra 0$, we obtain $$\Tor^{R}_{1}(R/I,\,S)\simeq\Ker(\delta) \mbox{\, and \,}  \Tor^{R}_{j+1}(R/I,\,S)\simeq\Tor^{R}_j(I,\,S)
\mbox{\, for all} \; j\ge 1,$$
where
$\delta:I\otimes_RS\to J$ is the multiplication map defined by
$x\otimes s\mapsto (x)\lambda s$ for $x\in I$ and $s\in S$. Clearly,
this implies that $(4)$ is  equivalent to $(3)$.

Now we show that $(1)$ and $(2)$ are equivalent.

According to Lemma \ref{epi00'}(1) and the fact that $\lambda$ is a
ring epimorphism, it follows that $\widetilde{\lambda}$ is a ring epimorphism. By
assumption, $J'$ is a left ideal of $S$, and therefore
$S\otimes_R(R/I)\simeq S/(S\cdot I)=S/J'$. Thanks to the general
result proved in the last part of the proof of \cite[Lemma 5.6]{xc3}, we
see
$$\Tor_i^{R/I}(S/J',\,W)\simeq \Tor_i^{R/I}(S\otimes_R(R/I),\,W)=0$$
for all $i\geq 1$ and all $S/J$-modules $W$. It then follows that
$\Tor_i^{R/I}(S/J',\,S/J)=0$ for all $i\geq 1$. Consider the short
exact sequence of right $R/I$-modules:
$$\quad 0\lra J/J'\lra S/J'\lra S/J\lra 0.$$
If we apply the functor $-\otimes_{R/I}(S/J)$ to this sequence, then
$\Tor^{R/I}_i(J/J',\,S/J)\simeq
\Tor^{R/I}_{i+1}(S/J,\,S/J)$ for all $i\geq 1$ and the
connecting homomorphism $\Tor^{R/I}_1(S/J,\,S/J)\to
(J/J')\otimes_{R/I}(S/J)$ is injective.

Clearly, if $\Tor^{R/I}_1(S/J,\,S/J)=0$, then
$\Tor^{R/I}_{j}(S/J,\,S/J)=0$ for all $j\geq 1$ if and only if
$\Tor^{R/I}_i(J/J',\,S/J)=0$ for all $i\geq 1$. This will imply that
$(1)$ and $(2)$ are
equivalent. So it is enough to demonstrate that $\Tor^{R/I}_1(S/J,\,S/J)=0$ always holds under
the assumptions of Corollary \ref{add}. However, this is true if we can show $(J/J')\otimes_{R/I}(S/J)=0$.

In fact,  if
$C\to D$ is a ring epimorphism, then $D\otimes_CX\simeq X$ as
$D$-modules for any $D$-module $X$, and $Y\otimes_CD\simeq Y$ as
right $D$-modules for any right $D$-module $Y$. This fact, together
with properties of ring epimorphisms, implies the following
isomorphisms:
$$(J/J')\otimes_{R/I}(S/J)\simeq (J/J')\otimes_{R}(S/J)\simeq
(J/J')\otimes_{R}\big(S\otimes_R(S/J)\big)\simeq\big((J/J')\otimes_{R}S\big)\otimes_R(S/J).$$
Since $SJ'=J'$ and $JJ'\subseteq J'$, we deduce
$\big((J/J')\otimes_{R}S\big)J'=0$. This means that
$(J/J')\otimes_{R}S$ is a right $S/J$-module. Clearly, the composite
of the two ring epimorphisms $R\to S$ and $S\to S/J$ is again a ring
epimorphism. It follows that
$\big((J/J')\otimes_{R}S\big)\otimes_R(S/J)\simeq (J/J')\otimes_R S$
as right $S/J$-modules.

In the following, we shall show $(J/J')\otimes_R S=0$.
Actually, applying the functor $-\otimes_RS$ to the exact sequence
$$0\lra J'\lraf{\alpha} J\lra J/J'\lra 0$$ of right $R$-modules, we
get an exact sequence
$$J'\otimes_RS\lraf{\alpha\otimes_RS} J\otimes_RS\lra
(J/J')\otimes_RS\lra 0$$ of right $S$-modules. Since $J$ is a right
$S$-module and $\lambda:R\to S$ is a ring epimorphism, the
multiplication map $\psi:J\otimes_RS\to J$, defined by $x\otimes
s\mapsto xs$ for $x\in J$ and $s\in S$, is an isomorphism. Note that
the map $(\alpha\otimes_RS)\psi: J'\otimes_RS\to J$ is surjective.
This yields that  $\alpha\otimes_RS$ is surjective and
$(J/J')\otimes_R S=0$. Hence $\Tor^{R/I}_1(S/J,\,S/J)=0$.
This finishes the proof. $\square$

\medskip
A special case of Lemma \ref{add} appears in trivial extensions. Let $\lambda: R\ra S$ be a homomorphism of rings and $M$ be an $S$-$S$-bimodule. Then $\lambda$ is homological if and only if $\widetilde{\lambda}: R\ltimes M\ra S\ltimes M$ is homological. The necessity of this condition follows from \cite[Theorem 1.1(1)]{xc3} and the proof of Corollary \ref{mod1b}. The sufficiency can be seen from Lemma \ref{add}.

Applying Theorem \ref{hom-dim} to the exact pair $(\lambda,\mu)$ in the proof of Lemma \ref{add}, we obtain the following estimations on finitistic dimensions, which can be applied to a class of examples of Milnor squares.

\begin{Koro} Let $\lambda: R\to S$ be a homological ring epimorphism. Suppose
that $I$ is an ideal of $R$ such that the image $J'$ of $I$ under
$\lambda$ is a left ideal in $S$ and that the restriction of
$\lambda$ to $I$ is injective. Let  $J$ be the ideal of $S$
generated by $J'$. Suppose that one of the conditions \emph{(1)-(4)} in Lemma \ref{add} holds. Then

$(1)$ $\fd(R)\leq \fd(S)+\fd(R/I)+\max\{1,\fld((R/I)_R)\}+1.$

$(2)$ If $_RS\in\Pf{R}$, then

\quad $(a)$ $\fd(S)\leq \fd(R)$ and $\fd(S/J)\leq \fd(R/I)$.

\quad $(b)$ $\fd(B)\leq \fd(R)+\fd(S/J)+\max\{1,\pd(_RS)\}+3$,
where $B:=\left(\begin{array}{lc} S & S/J'\\
0 & R/I\end{array}\right)$.
\end{Koro}

\medskip
{\bf Acknowledgement.} The both authors would like to express heartfelt gratitude to  Steffen Koenig for his kind help with the presentation of the results in this paper.  The research work of the corresponding author CCX is partially supported by NSF (KZ201410028033).

{\footnotesize
\bigskip Hongxing Chen,

School of Mathematical Sciences, BCMIIS, Capital Normal University, 100048
Beijing, People's Republic of  China

{\tt Email: chx19830818@163.com}

\bigskip
Changchang Xi,

School of Mathematical Sciences, BCMIIS, Capital Normal University, 100048
Beijing, People's Republic of  China

{\tt Email: xicc@cnu.edu.cn}}

\medskip
First version: June 26, 2013; Revised: April 24, 2014.


\medskip

{\footnotesize
\begin{thebibliography}{99}

\bibitem{AHKL}{{\sc L. angeleri-Huegel, S. Koenig} and {\sc Q. Liu,} Recollements and tilting objects, {\it J. Pure Appl. Algebra} {\bf 215} (2011) 420-438.}

\bibitem{AKLY}{{\sc L. Angeleri H\"ugel, S. K\"onig} , {\sc Q. Liu} and {\sc D. Yang},
Derived simple algebras and restrictions of recollements of derived module categories, arXiv:1310.3479v1, 2013.}

\bibitem{ars}{{\sc M. Auslander, I. Reiten} and {\sc S. Smalo},
\emph{Representation theory of Artin algebras}. Cambridge Studies in Advanced Mathematics 36.
Cambridge University Press, Cambridge, 1995.  }

\bibitem{bass}{{\sc H. Bass}, Finitistic dimension and a homological generalization of semiprimary rings,
\emph{Trans. Amer. Math. Soc.} \textbf{95} (1960) 466-488.}

\bibitem{BR}{{\sc A. Beligiannis} and {\sc I. Reiten},
Homological and homotopical aspects of torsion theories, {\it Mem.
Amer. Math. Soc}. {\bf 188} (2007), no. 883, 1-207.}

\bibitem{BBD}{{\sc A. A. Beilinson, J. Bernstein} and {\sc P.
Deligne},  Faisceaux pervers, \emph{Asterisque} \textbf{100} (1982)
5-171.}

\bibitem{xc1}{{\sc H. X. Chen} and {\sc C. C. Xi}, Good tilting
modules and recollements of derived module categories, \emph{Proc.
Lond. Math. Soc.}  104 (2012) 959-996.}

\bibitem{xc3}{{\sc H. X. Chen} and {\sc C. C. Xi}, Recollements of derived categories I: Exact
contexts, Preprint is available at:
http://math.bnu.edu.cn/$^{\sim}$ccxi/2012Test/Engindex.html, see also arXiv:1203.5168v2,
2012.}

\bibitem{xc5}{{\sc H. X. Chen} and {\sc C. C. Xi}, Recollements of derived categories II: Algebraic
$K$-theory, Primary version is available at:
 http://math.bnu.edu.cn/~ccxi/2012Test/Engindex.html, 2013. See also arXiv:1212.1879, 2012.}

 \bibitem{xc4}{{\sc H. X. Chen} and {\sc C. C. Xi}, Higher algebraic $K$-theory of ring epimorphisms, Preprint, 2013.
 Primary version is available at: http://math.bnu.edu.cn/$^{\sim}$ccxi/2012Test/Engindex.html, 2012.}


\bibitem{cohn}{{\sc P. M. Cohn}, On the free product of associative rings, \emph{Math. Z.}\textbf{ 71} (1959) 380-398.}

\bibitem{cohnbook1}{{\sc P. M. Cohn},  \emph{Free rings and their relations}, Academic Press, 1971.}

\bibitem{CPS}{{\sc E. Cline, B. Parshall} and {\sc L. Scott},
Derived categories and Morita theory, \emph{J. Algebra}
\textbf{104} (1986) 397-409.}

\bibitem{FGR}{{\sc R. Fossum, P. Griffith} and {\sc I. Reiten},
\emph {Trivial extensions of abelian categories with applications to ring theory}, Lecture Notes in Math. {\bf 456}, Springer, Berlin, 1975.}

\bibitem{Happel1}{{\sc D. Happel}, {\it Triangulated categories in the representation theory of finite dimensional algebras}.
London Math. Soc. Lecture Note Series \textbf{119}, Cambridge
University Press, 1988.}


\bibitem{Happel}{{\sc D. Happel}, Reduction techniques for homological conjectures, \emph{Tsukuba
J. Math.} \textbf{17} (1993), no. 1, 115-130.}

\bibitem{hx2}{{\sc W. Hu } and {\sc C. C. Xi},
$\cal D$-split sequences and derived equivalences, \emph{Adv. Math.}
\textbf{227} (2011) 292-318.}

\bibitem{koenig}{{\sc S. Koenig}, Tilting complexes, perpendicular categories and recollements of derived module categories of rings,
\emph{J. Pure Appl. Algebra} \textbf{73} (1991) 211-232.}

\bibitem{Milnor}{{\sc J. Milnor}, \emph{Introduction to algebraic $K$-theory}.
Annals of Mathematics Studies, vol. \textbf{72}.
Princeton University Press, Princeton (1971)}

\bibitem{nr}{{\sc A. Neeman} and {\sc A. Ranicki}, Noncommutative localization in algebraic $K$-theory I,
\emph{Geometry and Topology} \textbf{8} (2004) 1385-1425. }

\bibitem{NS2}{{\sc P. Nicol\'as} and {\sc M. Saor\'in},
Lifting and restricting recollement data, {\it Appl. Categor. Struct.} {\bf 19}(3) (2011) 557-596.}

\bibitem{XP}{{\sc S. Y. Pan} and {\sc C. C. Xi,} Finiteness of finitistic dimension is invariant
under derived equivalences, {\it J. Algebra}
\textbf{322}(2009) 21-24.}

\bibitem{x2}{{\sc C. C. Xi}, On the finitistic dimension conjecture II: Related to finite global dimension,
{\it Adv. Math.} {\bf 201} (2006) 116-142.}

\bibitem{xx}{{\sc C. C. Xi} and {\sc D. M. Xu}, The finitistic dimension conjecture and relatively projective modules,
{\it Commun. Contemp. Math.} \textbf{15} (2) (2013) 1350004(27 pages). DOI: 10.1142/S0219199713500041.}

\bibitem{DMX}{{\sc D. M. Xu}, Homological dimensions and strongly
idempotent ideals, {\it J. Algebra} (to appear). Prerint, arXiv:1303.1306v1, 2013.}

\bibitem{ZH}{{\sc B. Zimmermann-Huisgen}, Homological domino effects and the first finitistic dimension conjecture, {\it Invent. Math.} {\bf 108} (1992) 369-383.}

\end{thebibliography}}
\end{document}